\documentclass[11pt,a4,onesided]{amsart}
 \usepackage{epsfig}
\usepackage{amssymb}
\usepackage{color}
\usepackage{mathrsfs}
\usepackage{hyperref}

\newcommand\R{{\mathbb R}}

\def\AA{{\mathcal A}}
\def\BB{{\mathcal B}}

\def\DD{{\mathcal D}}

\def\FF{{\mathcal F}}
\def\GG{{\mathcal G}}
\def\HH{{\mathcal H}}

\def\LL{{\mathcal L}}
\def\MM{{\mathcal M}}

\def\SS{{\mathcal S}}
\def\TT{{\mathcal T}}

\newcommand{\sk}{\smallskip}

\let\oldmarginpar\marginpar
\renewcommand\marginpar[1]{\-\oldmarginpar[\raggedleft\footnotesize #1]%
{\raggedright\footnotesize #1}}

\newtheorem{theo}{Theorem}

\newtheorem{lem}[theo]{Lemma}

\theoremstyle{definition}

\theoremstyle{remark}
\newtheorem{rem}[theo]{Remark}

\newcommand{\beqn}{\begin{equation}}
\newcommand{\eeqn}{\end{equation}}
\newcommand{\bear}{\begin{eqnarray}}
\newcommand{\eear}{\end{eqnarray}}
\newcommand{\bean}{\begin{eqnarray*}}
\newcommand{\eean}{\end{eqnarray*}}

\begin{document}

\title[The KFP Equation with weak confinement]{The Kinetic Fokker-Planck Equation with Weak Confinement Force}

\author{Chuqi Cao }

\maketitle

\begin{center} {\bf 
Version of \today}
\end{center}

\begin{abstract} 

We consider the kinetic Fokker-Planck equation with weak confinement force. We proved some (polynomial and sub-exponential) rate of convergence to the equilibrium (depending on the space to which the initial datum belongs). Our results generalized the result in \cite {DMS, DMS2, V, HN, H, HN2, BCG, MM} to weak confinement case.

  \end{abstract}

\bigskip

Asymptotic behavior of solutions
 [see also 45C05, 45K05,  35410],
35Q84 
Fokker-Planck equations.
\bigskip


\bigskip

\tableofcontents

\section{Introduction}
\label{sec1}
\setcounter{equation}{0}
\setcounter{theo}{0}

In this paper, we consider the weak hypocoercivity issue for the kinetic Fokker-Planck (KFP for short) equation 
\sk
\beqn\label{1}
\partial_t f =\LL f := -v \cdot \nabla_x f +\nabla_x V(x) \cdot \nabla_v f +\Delta_v f + \hbox{div}_v(v f),
\eeqn
for a density function $f = f(t, x, v)$, with $t \ge 0 ,\ x \in \R^d,\ v \in \R^d$.
The evolution equation is complemented with an initial datum
\beqn
\nonumber
f(0, \cdot) = f_0 \ \ on \ \R^{2d} .
\eeqn
We make the fundamental assumption on the  confinement potential $V$
\beqn
\nonumber
V(x)= \langle x \rangle^{\gamma} ,\quad \gamma \in (0, 1) ,
\eeqn
where $ \langle x \rangle^2 := 1 + |x|^2 $.

Let us make some elementary but fundamental observations. First, the equation is mass conservative, 
that is
\beqn
\nonumber
\MM (f_0) =\MM( f(t,\cdot )),
\eeqn
where we define the mass of $f$ by
\beqn
\nonumber
\MM(f) =\int_{R^d \times R^d} f dx dv.
\eeqn

\noindent Next, we observe that 
\beqn\label{A1}
\quad G=Z^{-1}e^{-W}, \quad W= \frac {v^2} {2} +V(x), \quad Z \in \R_{+}
\eeqn
is a positive normalized steady state of the KFP model, precisely 
\beqn
\nonumber
LG=0, \quad G>0, \quad \MM(G)=1, 
\eeqn
by choosing the normalizing constant $Z>0$ appropriately. Finally we observe that, contrary to the case $\gamma \ge 1$, a Poincar\'{e} inequality of the type
\beqn
\nonumber
\exists c > 0, \quad \int_{\R^d} |f(x)|^2 \exp(-V(x)) dx  \le c \int_{\R^d} |\nabla f(x)|^2 \exp(-V(x)) dx ,
\eeqn
for any smooth function $f:  \R^d \to \R$ such that 
\beqn
\nonumber
\int_{\R^d}  f(x) \exp(-V(x)) dx =0 ,
\eeqn
does not hold. Only a weaker version of this inequality remains true (see \cite{RW}, or below Section 2). In particular, there is no spectral gap for the associated operator $\LL$, nor is there an exponential trend to the equilibrium for the associated semigroup.   

\sk

For a given weight function $m$, we will denote $L^p(m) = \{ f | fm \in L^p \}$ the associated Lebesgue space and $\Vert f \Vert_{L^p(m)} = \Vert f m \Vert_{L^p}$ the associated norm.

The notation $A \lesssim B$ means $A \le C  B $ for some constant $C>0$.
\smallskip

With these notations, we can introduce the main result of this paper.

\begin{theo}\label{T11}

(1)  For any initial datum $f_0 \in L^p(G^{-(\frac {p-1} {p} +\epsilon)})$, $p \in [1,\infty)$, $\epsilon > 0$ small, the associated solution $f(t, \cdot)$ of the kinetic Fokker-Planck equation (\ref{1}) satisfies
\beqn
\nonumber
\Vert f(t , \cdot)	- \MM(f_0) G\Vert_{L^p(G^{-\frac {p-1} {p}})}  \lesssim e^{- C t^{b} } \Vert f_0 -\MM(f_0) G\Vert_{L^p (G^{-(\frac {p-1} {p} +\epsilon) } ) }  ,
\eeqn
for any $b \in (0, \frac {\gamma} {2-\gamma})$ and some constant $C>0$.\\
(2) For any initial datum $f_0 \in L^1(m)$, $m = H^k$, $H =x^2+v^2$, $1 \le k $, the associated solution $f(t, \cdot)$ of the kinetic Fokker-Planck equation (\ref{1}) satisfies
\bear
\nonumber
\Vert f(t,\cdot)-\MM(f_0)G \Vert_{L^1} \lesssim  (1+t)^{-a} \Vert f_0 - \MM(f_0)G \Vert_{L^1(m)},
\eear
for any $0 < a <\frac {k} {1-\frac \gamma 2}$. The constants in the estimates do not depend on $f_0$,  but rely on $\gamma, d, \epsilon,\theta, p, k$.
\end{theo}

\begin{rem}

Theorem \ref{T11} is also true when $V(x)$ behaves like $ \langle x \rangle^{\gamma} $, that is for any $V(x)$ satisfying

\beqn
\nonumber
C_1 \langle x \rangle^{\gamma} \le V(x) \le C_2\langle x \rangle^{\gamma} , \ \ \forall x \in \R^d ,
\eeqn

\beqn
\nonumber
C_3|x| \langle x \rangle^{\gamma-1} \le x \cdot \nabla_x V(x) \le C_4|x|\langle x \rangle^{\gamma-1}, \ \ \forall x \in B_R^c ,
\eeqn
and 
\beqn
\nonumber
|D^2_x V(x)| \le C_5 \langle x \rangle^{\gamma-2}, \ \ \forall x \in \R^d ,
\eeqn
for some constant $C_i >0$, $R>0$.
\end{rem}

\begin{rem}

There are many classical results on the case $ \gamma \ge 1 $. In this case there is an exponentially decay, and we refer the interested readers to \cite{V, DMS, DMS2, HN, HN2, H, BCG}.

\end{rem}

\begin{rem}
There are already some convergence results for the weak confinement case proved by probability method on some particular $L^1$ or $L^2$ spaces in \cite{BCG} and \cite{DFG}, this paper extend the result to $L^p$ spaces and more larger spaces.
\end{rem}
\sk

Let us briefly explain the main ideas behind our method of proof. 

\sk

We first introduce four spaces $E_1= L^2(G^{-1/2})$, $E_2=L^2(G^{-1/2}e^{\epsilon_1V(x)})$, $E_3 = L^2(G^{-(1+\epsilon_2)/2})$ and $E_0 =L^2(G^{-1/2} \langle x \rangle^{\gamma-1} )$, with 
 $\epsilon_1>0$ and $\epsilon_2>0$ small such that $E_3 \subset E_2 \subset E_1 \subset E_0 \subset L^2$. Thus $E_1$ is an interpolation space between $E_0$ and $E_2$. We first use a hypocorecivity argument as in \cite{DMS,DMS2} to prove that, for any $f_0 \in E_3$, the solution to the KFP equation (\ref{1}) satisfies 
\beqn
\nonumber
\frac {d} {dt} \Vert f(t) \Vert_{E_1} \le - \lambda \Vert f(t) \Vert_{E_0},
\eeqn
for some constant $\lambda >0$. We use this and the Duhamel formula to prove
\beqn
\nonumber
\Vert f(t) \Vert_{E_2} \lesssim \Vert f_0 \Vert_{E_3}.
\eeqn
Combining  the two inequalities and using a interpolation argument as in \cite{KM}, we get
\beqn\label{7}
\Vert f(t) \Vert_{E_1} \lesssim e^{-a t^b} \Vert f_0 \Vert_{E_3},
\eeqn
for some $a>0, b \in (0,1)$.

We then generalize the decay estimate to a wider class of Banach spaces by adapting the extension theory introduced in \cite{M2} and developed in \cite{MM, GMM}.	
For any operator $\LL$, denote $S_\LL(t)$ the associated semigroup. We introduce a splitting $\LL=\AA+\BB$, where $\AA$ is an appropriately defined bounded operator so that $\BB$ becomes a dissipative operator. By proving  some regularization estimate in $S_\BB$ in $L^p$
\beqn
\nonumber
\Vert S_\BB(t) \Vert_{L^p(m_1) \to L^2(m_2)} \lesssim t^{-\alpha},\quad \forall t \in [0,\eta],
\eeqn
for some weight function $m_1$, $m_2$ and some $\alpha, \eta>0$, and using the iterated Duhamel's formula
\bear\label{10}
S_\LL =S_\BB + \sum_{l=1}^{n-1}(S_\BB)*( \AA S_\BB )^{(*l)}  + S_\LL *(\AA S_\BB(t))^{*n},
\eear
we deduce the $L^p$ convergence on $S_\LL$, where the convolution of two semigroups $S_\AA(t)$ $S_\BB(t)$is defined by
\beqn
\nonumber
(S_\AA*S_\BB)(t) =\int_0^t S_\AA(s) S_\BB(t-s) ds.
\eeqn

\medskip
Let us end the introduction by describing the plan of the paper. In Section 2, we will develop a hypocoercivity argument to prove a weighted $L^2$ estimate for the KFP model.  In section 3, we introduce a splitting $\LL =\AA+\BB$ and using the $L^2$ estimate, we prove a $L^2$ convergence. 
In Section 4 we present the proof of a regularization estimate on $S_\BB$ from $L^2$ to $L^p$. In Section 5 we prove some $L^1$ estimate on the semigroup $S_\BB$. Finally in Section 6 we use the above regularization estimate to conclude the  $L^p$ convergence for KFP equation.

\medskip
{\bf Acknowledgment. }Â Â 
The author thanks to S. Mischler for furitful discussions on the full work of the paper. This work was supported by grants from R\'egion Ile-de-France the DIM program.

\bigskip
\section{$L^2$ framework: Dirichlet form and rate of convergence estimate}
\label{sec2}
\setcounter{equation}{0}
\setcounter{theo}{0}

For later discussion, we introduce some notations for the whole paper.

\sk 

\noindent We split the KFP operator  as $$ \LL = \TT + \SS, $$ where $\TT$ stands for the transport part
\beqn
\nonumber
\TT f= -v \cdot \nabla_x f +\nabla_x V(x) \cdot \nabla_v f,
\eeqn
and $\SS$ stands for the collision part
\beqn
\nonumber
\SS f= \Delta_v f + div_v(v f).
\eeqn
We will denote the cut-off function $\chi$ such that  $\chi(x, v) \in [0, 1]$, $\chi(x, v) \in C^\infty $, $\chi(x, v) =1$ when $x^2+v^2 \le 1$ , $\chi(x, v) =0$ when $x^2+v^2 \ge 2$, and then denote $\chi_R = \chi(x/R, v/R)$.

\sk
\noindent We may also define another splitting of the KFP operator $\LL$ by
\beqn\label{21}
\LL =\AA+\BB, \quad \AA=K\chi_R(x, v).
\eeqn
with $K, R >0$ to be chosen later.

\noindent We use $\int f$ in place of $\int_{\R^d \times R^d} f dx dv$ for short, similarly $\int f dx$ means $\int_{\R^d } f dx $ , $\int f dv $  means $\int_{\R^d } f dv$. $B_{|x| \le \rho}$ is used to denote the ball such that $\{ x\in \R^d | |x| \le \rho\} $, similarly $B_{\rho}$ means the ball such that $\{ x, v \in \R^d | |x|^2 +v^2 \le \rho\} $.

For $V(x)= \langle x \rangle^\gamma, 0 < \gamma < 1$, we also denote $\langle \nabla V \rangle$ for $\langle x \rangle^{\gamma-1}$, and $\langle \nabla V \rangle^{-1}$ for $\langle x \rangle^{1-\gamma}$.

\sk

With these notations we introduce the Dirichlet form adapted to our problem. We define the 0 order and first order moments
\beqn
\nonumber
\rho_f=\rho[f]=\int f dv, \quad  j_f = j[f] = \int v  f dv,
\eeqn
then we define a projection operator $\pi$ by
\beqn
\nonumber
 \pi f =M \rho_f, \quad M=C e^{-v^2/2}, \quad \int M dv=1,
\eeqn
and the complement of $\pi$ by
\beqn
\nonumber
  \pi^{\bot} = 1 - \pi, \quad f^{\bot} = \pi^{\bot} f.
\eeqn
We define an elliptic operator $\Delta_V$ and its dual $\Delta_V^*$ by
\beqn
\nonumber
\Delta_V u :=div_x(\nabla_x u + \nabla_x V u),\quad \Delta_V^* u =\Delta_x u- \nabla_x V \cdot \nabla_x u,
\eeqn
let $u=(\Delta_V^*)^{-1}\xi$ be the solution to the above elliptic equation 
\beqn
\nonumber
\Delta_V^* u =\xi \ on \ \R^d,
\eeqn
note that u can differ by a constant, we also requires that
\beqn
\nonumber
\int u e^{-V} \langle \nabla V \rangle^{-2} dx=0,
\eeqn
using these notations, define a scalar product by
\bear
\nonumber
((f, g))&:= &(f, g)_{\HH} +\epsilon (\Delta_V^{-1} \nabla_x j_f, (\rho_g e^V \langle \nabla V \rangle^2 ))_{L^2} 
\\ \nonumber
&& + \epsilon ((\rho_f e^V \langle \nabla V \rangle^2), \Delta_V^{-1} \nabla_x j_g)_{L^2}
\\ \nonumber
&=&(f, g)_{\HH} + \epsilon (j_f, \nabla_x(\Delta_V^*)^{-1}(\rho_ge^V \langle \nabla V \rangle^2 ))_{L^2} 
\\ \nonumber
&&+\epsilon( (\nabla_x(\Delta_V^*)^{-1}(\rho_f e^V \langle \nabla V \rangle^2 ), j_g)_{L^2},
\eear
for some $\epsilon > 0$ to be specified later. \\
We then define the Dirichlet form 
\bear
\nonumber
D[f] &:= &((-\LL f, f))
\\ \nonumber
&=& (-\LL f, f)_{\HH} +\epsilon (\Delta_V^{-1} \nabla_x j[-\LL f], (\rho_f e^V \langle \nabla V \rangle^2 ))_{L^2}
\\ \nonumber
&&+ \epsilon ((\rho[-\LL f] e^V \langle \nabla V \rangle^2), \Delta_V^{-1} \nabla_x j_f)_{L^2}.
\eear
Finally we define $\HH=L^2(G^{-1/ 2 })$, $\HH_1=L^2(G^{-1/ 2 }\langle \nabla V \rangle)$ and 
\beqn
\nonumber
\HH_0 = \{h \in \HH , \int f dx dv =0  \}
\eeqn
where we recall that $G$ has been introduced in (\ref{A1}). With these notations we can come to our first theorem.
\begin{theo}\label{T21}
 There exists $\epsilon > 0$ small enough, such that on $\HH_0$ the norm $((f,  f))^{\frac 1 2} $ defined above is equivalent to the norm of $\HH$, moreover there exist $\lambda > 0$, such that
$$D[f]  \ge \lambda \Vert f \Vert^2_{\HH_1}, \ \  \ \forall f \in \HH_0.$$
As a consequence, for any $f_0 \in \HH_0$, we have
\beqn\label{22}
\frac d {dt} ((f, f)) \le  - C \int f^2 G^{-1} \langle x \rangle^{2(\gamma - 1)},
\eeqn
for some constant $C>0$. In particular for any $f_0 \in \HH_0$, we have
\beqn\label{23}
\Vert f(t , \cdot) \Vert_{L^2(G^{-\frac 1 2})} \le C \Vert f_0 \Vert_{L^2(G^{-\frac 1 2})},
\eeqn
for some constant $C > 0$.
\end{theo}

\begin{rem}
In $\HH_0$ we have
\beqn
\nonumber
\int  \rho_f e^V \langle \nabla V \rangle^2  e^{-V} \langle \nabla V \rangle^{-2}  dx= \int \rho_f dx = \int f dx dv = 0,
\eeqn
so the term $(\Delta_V^*)^{-1}(\rho_ge^V \langle \nabla V \rangle^2 )$  is well defined in $\HH_0$.
\end{rem}

\begin{rem}
(1) By little modifying the method in Villani's paper \cite{V}, a $H^1$ version of our theorem can be established.\\
(2) Our statement is a generalization of \cite{DMS, DMS2}.
\end{rem}

Before proving the theorem, we need some lemmas. 

We say that $W$ satisfies a local Poincar\'e inequality on a bounded open set $\Omega$ if there exist some constant $\kappa_\Omega > 0$ such that:
\beqn
\nonumber
\int_\Omega h^2 W \le k_\Omega \int_\Omega |\nabla h |^2 W + \frac 1 {W(\Omega)}(\int_\Omega h W)^2,
\eeqn
for any nice function $h : \R^d \to \R$ and where we denote $W(\Omega):= \langle W 1_\Omega \rangle$.
\begin{lem}\label{L21}
Under the assumption $W, W^{-1} \in L_{loc}^\infty(\R^d) $, the function $W$ satisfies the local Poincar\'e inequality for any ball $\Omega \in \R^d$.
\end{lem}
\noindent For the proof of Lemma \ref{L21} we refer to \cite{M3} Lemma 2.3.

\begin{lem}\label{L22}(weak Poincar\'{e} inequality)
There exist a constant $\lambda > 0$ such that 
\beqn
\nonumber
\Vert u \Vert_{L^2(\langle \nabla V \rangle e^{-V/2}  )} \le \lambda \Vert \nabla u \Vert_{L^2( e^{-V/2} )}
\eeqn
for any $u \in \DD(\R^d)$ such that 
\beqn
\nonumber
\int_{\R^d}  u e^{-V}\langle \nabla V \rangle^{-2}  dx =0
\eeqn
\end{lem}

\noindent {\sl Proof of Lemma \ref{L22}.} We prove for any $h \in \DD(\R^d)$ such that
\beqn\label{24}
\int_{\R^d} h e^{-V} \langle \nabla V \rangle^{-2}=0,
\eeqn
we have
\beqn
\nonumber
\int_{\R^d} |\nabla h |^2  e^{-V} \ge \lambda \int_{\R^d} h^2 e^{-V} \langle x\rangle^{2(\gamma-1)},
\eeqn
for some $\lambda >0$. Taking $g=he^{-1/2V}$, we have $\nabla g =\nabla h e^{-\frac 1 2 V} -\frac 1 2 \nabla V h e^{-\frac 1 2 V}$, so that
\bear
\nonumber
0 \le\int |\nabla g|^2 &=& \int |\nabla h |^2e^{-V} +\int h^2 \frac 1 4 |\nabla V|^2 e^{-V} -\int \frac 1 2\nabla(h^2) \cdot \nabla V e^{-V}
\\ \nonumber
&=&\int |\nabla h |^2e^{-V} +\int h^2 (\frac 1 2 \Delta V - \frac 1 4 |\nabla V|^2 )e^{-V} .
\eear
We deduce for some $K, R_0 > 0$
\bear
\nonumber
\int |\nabla h |^2e^{-V} \ge \int \frac 1 8 h^2 \langle \nabla V \rangle^2 e^{-V} - K \int_{B_{R_0}} h^2 e^{-V} \langle \nabla V \rangle^{-2}.
\eear
Defining
\beqn
\nonumber
\epsilon_R :=\int_{B_R^c} e^{-V} \langle \nabla V \rangle^{-6} , \quad Z_R :=\int_{B_R} e^{-V} \langle \nabla V \rangle^{-2} ,
\eeqn
and using (\ref{24}), we get
\bear
\nonumber
(\int_{B_R} he^{-V} \langle \nabla V \rangle^{-2} )^2&=&(\int_{B_R^c} he^{-V} \langle \nabla V \rangle^{-2} )^2
\\ \nonumber
&\le&\int_{B_R^c} h^2e^{-V} \langle \nabla V \rangle^{2} \int_{B_R^c} e^{-V} \langle \nabla V \rangle^{-6} 
\\ \nonumber
&\le& \epsilon_R \int_{B_R^c} h^2e^{-V} \langle \nabla V \rangle^{2}.
\eear
Using the local Poincar\'e inequality in Lemma \ref{L21}, we deduce 
\bear
\nonumber
\int_{B_R} h^2e^{-V} \langle \nabla V \rangle^{-2} &\le& C_R \int_{B_R} |\nabla h|^2e^{-V} \langle \nabla V \rangle^{-2} +\frac 1 {Z_R} (\int_{B_R} he^{-V} \langle \nabla V \rangle^{-2} )^2
\\ \nonumber
&\le& C_R^{'} \int_{B_R} |\nabla h|^2e^{-V}  +\frac {\epsilon_R} {Z_R} \int_{B_R} h^2e^{-V} \langle \nabla V \rangle^{2} .
\eear
Putting all the inequalities together and taking $R>R_0$, we finally get
\bear
\nonumber
\int h^2e^{-V} \langle \nabla V \rangle^{2} &\le& 8\int |\nabla h |^2e^{-V} +8 K \int_{B_{R_0}} h^2 e^{-V} \langle \nabla V \rangle^{-2}
\\ \nonumber
&\le& 8(1+ KC_R^{'} )\int |\nabla h|^2e^{-V}  +\frac {8 K \epsilon_R} {Z_R} \int_{B_R} h^2e^{-V} \langle \nabla V \rangle^{2} ,
\eear
and we conclude by taking $R$ large such that: $\frac {8 K \epsilon_R} {Z_R} \le \frac 1 2  $.
\qed

\smallskip

\begin{lem}\label{L23}

(Elliptic Estimate) For any $\xi_1 \in L^2( \langle \nabla V \rangle^{-1}e^{-V/2})$ and $\xi_2 \in L^2(e^{-V/2})$, the solution $u \in L^2(e^{-V/2})$ to the elliptic equation

\beqn\label{E21}
-\Delta_V^* u = \xi_1 + \nabla \xi_2, \quad \int u e^{-V} \langle \nabla V \rangle^{-2}dx =0,
\eeqn

satisfies

\beqn\label{E22}
\Vert u \Vert_{L^2(\langle \nabla V \rangle e^{-V/2} )}+ \Vert \nabla u \Vert_{L^2(e^{-V/2})} \lesssim \Vert \xi_1 \Vert_{L^2(\langle \nabla V \rangle^{-1} e^{-V/2})} + \Vert \xi_2 \Vert_{L^2(e^{-V/2})}.
\eeqn

Similarly for any $\xi \in L^2(e^{-V/2})$, the solution $u \in L(e^{-V/2}) $ to the elliptic problem

\beqn
\nonumber
-\Delta_V^* u = \xi, \quad \int u e^{-V} \langle \nabla V \rangle^{-2} =0,
\eeqn

satisfies

\beqn\label{E23}
\Vert u \Vert_{L^2(\langle \nabla V \rangle^2 e^{-V/2} )}+ \Vert \nabla u \Vert_{L^2(\langle \nabla V \rangle e^{-V/2})} + \Vert D^2 u \Vert_{L^2(e^{-V/2})} \lesssim  \Vert \xi \Vert_{L^2(e^{-V/2 }\langle \nabla V \rangle^{-1} ) }.
\eeqn

\end{lem}

\noindent {\sl Proof of Lemma \ref{L23}.}
Multiply (\ref{E21}) by $u e^{-V}$ and observes that
\beqn\label{E24}
e^{V} \hbox{div}_x[e^{-V}\nabla_x u] = \Delta_x u -\nabla_x V \cdot \nabla_x u= \Delta_V^* u,
\eeqn
we have after integration
\beqn
\nonumber
- \int e^V \hbox{div}_x[e^{-V}\nabla_x u] u e^{-V} = \int (\xi_1 + \nabla \cdot \xi_2) u e^{-V}.
\eeqn
Performing one integration by parts, we deduce
\beqn
\nonumber
\int e^{-V} |\nabla_x u|^2 = \int( \xi_1 u +\xi_2 \cdot \nabla u -\xi_2 \cdot \nabla V u )e^{-V},
\eeqn
using Lemma \ref{L22} we obtain (\ref{E22}). In inequality (\ref{E23}), the first two terms are easily bounded by (\ref{E22}) and $\langle \nabla V \rangle \le 1$, we then only need to prove the bound for the third term. By integration by parts, we have
\bear
\nonumber
\int |D^2 u|^2e^{-V} &=& \sum\limits_{i, j=1}^d \int (\partial_{ij}^2 u)^2e^{-V} 
\\ \nonumber
&=& \sum\limits_{i, j =1}^d \int \partial_i u (\partial_{ij}^2 u \partial_j V -\partial_{ijj}^3 u  )e^{-V}
\\ \nonumber
&=& \sum\limits_{i, j=1}^d \int \partial_{jj} u \partial_i(\partial_i u e^{-V}) -\frac 1 2 \int (\partial_i u)^2 \partial_j(\partial_j V e^{-V})
\\ \nonumber
&=&  \int (\Delta u) (-\Delta_V^* u) e^{-V}  + \int |\nabla u|^2( |\nabla V|^2 -\Delta V )e^{-V}
\\ \nonumber  
&\lesssim&\Vert D^2 u\Vert_{L^2(e^{-V/2})} \Vert \xi \Vert_{L^2(e^{-V/2}) } + \Vert \langle \nabla V\rangle \nabla u \Vert_{L^2(e^{-V/2})},
\eear
where in the third equality we have used 
\bear
\nonumber
\int \partial_{ij}^2 u\partial_i u  \partial_j V e^{-V} &=& -\int \partial_i u \partial_j (\partial_i u \partial_j Ve^{-V})
\\ \nonumber
&=& -\int \partial_{ij}^2 u\partial_i u  \partial_j V e^{-V} -\int (\partial_i u)^2 \partial_j(\partial_j V e^{-V}),
\eear
which implies 
\bear
\nonumber
 \int \partial_{ij}^2 u\partial_i u  \partial_j V e^{-V} = -\frac 1 2 \int (\partial_i u)^2 \partial_j(\partial_j V e^{-V}),
\eear
and in the fourth equality we have used (\ref{E24}). That concludes the proof.
\qed

\smallskip

Now we turn to the proof of Theorem \ref{T21}.

\sk

\noindent {Proof of Theorem \ref{T21}.} First we prove the equivalence of the norms associated to $((\ , \ ))$ and $(\  ,\  )_\HH$. By Cauchy-Schwarz inequality and Lemma \ref{L21}, we have
\beqn
\nonumber
(j_f, \nabla_x(\Delta_V^*)^{-1}(\rho_ge^V \langle \nabla V \rangle^2 ))_{L^2} \le \Vert j_f  \Vert_{L^2(e^{V/2})} \Vert \rho_g e^V \langle \nabla V \rangle^2 \Vert_{L^2(\langle \nabla V \rangle^{-1} e^{-V/2})},
\eeqn
and obviously
\beqn
\nonumber
\Vert \rho_g e^V \langle \nabla V \rangle^2 \Vert_{L^2(\langle \nabla V \rangle^{-1} e^{-V/2})} = \Vert \rho_g \Vert_{L^2( \langle \nabla V \rangle e^{V/2})} \le \Vert \rho_g \Vert_{L^2( e^{V/2}) } \lesssim \Vert g \Vert_{\HH}.
\eeqn
Using the elementary observations
\beqn
\nonumber
|j_f |\lesssim \Vert f \Vert_{L^2(e^{v^2/4})} \ \ \ |\rho_f |\lesssim \Vert f \Vert_{L^2(e^{v^2/4})},
\eeqn
we deduce
\beqn
\nonumber
(j_f, \nabla_x(\Delta_V^*)^{-1}(\rho_ge^V \langle \nabla V \rangle^2 ))_{L^2} \lesssim \Vert f \Vert_{\HH}\Vert g \Vert_{\HH},
\eeqn
The third term in the definition of (( , ))  can be estimated in the same way and that ends the proof of  equivalence of norms. 
\qed

\sk

Now we prove the main estimate of the theorem. We split the Dirichlet term $D[f]$ into 3 parts
$$D[f]=T_1+\epsilon T_2 + \epsilon T_3,$$
with
\bear
\nonumber
T_1&:=&(\LL f , f)_{\HH}
\\ \nonumber
T_2 &:=&  (\Delta_V^{-1} \nabla_x j[-\LL f],  \rho_f)_{L^2(e^{V/2 }  \langle \nabla V \rangle )}
\\ \nonumber
T_3 &:=&( (\Delta_V)^{-1}\nabla_x j_f, \rho[-\LL f])_{L^2(e^{V/2} \langle \nabla V \rangle)} \ ,
\eear
and compute them separately.\\
For the $T_1$ term, using the classical Poincar\'e inequality, we have
\bear
\nonumber
T_1 &:=&(-\TT f+ \SS f, f)_{\HH}=(-\SS f, f)_{\HH}
\\ \nonumber
&= &-\int [\Delta_v f + div_v( v f)] f M^{-1} e^{V} = \int  |\nabla_v (f /M)|^2 M e^{V}
\\ \nonumber
&\ge&   k_p \int | f/M - \rho_f |^2M e^{V}  = k_p \Vert f - \rho_f M \Vert^2_{\HH} = k_p \Vert f^{\bot} \Vert_{\HH}^2,
\eear
for some $k_p>0$. We split the $T_2$ term as
\bear
\nonumber
T_2 &:=&  (\Delta_V^{-1} \nabla_x j[-\LL f],  \rho_f)_{L^2(e^{V/2 }  \langle \nabla V \rangle )}
\\ \nonumber
&=&(\Delta_V^{-1} \nabla_x j[-\TT \pi f],  \rho_f)_{L^2(e^{V/2}  \langle \nabla V \rangle  )}
\\ \nonumber
&&+  (\Delta_V^{-1} \nabla_x j[-\TT f^{\bot}],  \rho_f)_{L^2(e^{V/2}  \langle \nabla V \rangle )}
\\ \nonumber
&&+  (\Delta_V^{-1} \nabla_x j[-\SS f],  \rho_f)_{L^2(e^{V/2}  \langle \nabla V \rangle  )}
 \\ \nonumber
&:=& T_{2,1} + T_{2,2} +T_{2,3}  .
\eear
First observe
\beqn
\nonumber
\TT \pi f = - v \cdot \nabla_x \rho_f M - \nabla_x V  \cdot v \rho_f M = -e^{-V} M v \cdot \nabla_x (\rho_f/e^{-V}),
\eeqn
so that we have
\beqn
\nonumber
j[-\TT \pi f]=\langle v v_k M\rangle e^{-V}\partial_{x_k}(\rho_f/e^{-V})=e^{-V}\nabla_x(\rho_f/e^{-V}).
\eeqn
Next by (\ref{E24}), we have
\bear
\nonumber
T_{2,1} &=&(j[-\TT \pi f], \nabla (\Delta_V^*)^{-1}(\rho_f e^V \langle \nabla V \rangle^2 ))_{L^2}
\\ \nonumber
&=&(\rho_f, [e^Vdiv_x(e^{-V}\nabla)][(\Delta_V^*)^{-1}(\rho_f e^V \langle \nabla V \rangle^2 )] )_{L^2}
\\ \nonumber
&=& \Vert \rho_f e^{V/2} \langle \nabla V \rangle  \Vert^2_{L^2} = \Vert \pi f \Vert^2_{\HH_1}.
\eear
Using the notation $\eta_{1} = \langle v \otimes v f^{\bot} \rangle$ and $\eta_{2,\alpha \beta} = \langle v_\alpha \partial_{v_\beta} f^{\bot} \rangle$, and observing that
\beqn
\nonumber
|\eta_1|  \lesssim \Vert f^{\bot} \Vert_{L^2(e^{v^2/4})}    ,   |\eta_2|  \lesssim \Vert f^{\bot} \Vert_{L^2(e^{v^2/4})},
\eeqn
we  compute
\bear
\nonumber
T_{2,2} &=&(j[-\TT f^{\bot}], \nabla(\Delta_V^*)^{-1} (\rho_f e^V \langle \nabla V \rangle^2 ))_{L^2}
\\ \nonumber
&=& (D\eta_1+ \eta_2 \nabla V, \nabla(\Delta_V^*)^{-1} (\rho_f e^V \langle \nabla V \rangle^2 ))_{L^2}
\\ \nonumber
&=&  (\eta_1, D^2(\Delta_V^*)^{-1} (\rho_f e^V \langle \nabla V \rangle^2 ))_{L^2} + (\eta_2 ,\nabla V \nabla(\Delta_V^*)^{-1} (\rho_f e^V \langle \nabla V \rangle^2 ))_{L^2}
\\ \nonumber
&= &\Vert \eta_1 \Vert_{L^2(e^{V/2}  )}\Vert D^2(\Delta_V^*)^{-1} (\rho_f e^V\langle \nabla V \rangle^2) \Vert_{L^2(e^{-V/2})} 
\\ \nonumber
&&+ \Vert \eta_2 \Vert_{L^2(e^{V/2})}\Vert \nabla V \nabla(\Delta_V^*)^{-1} (\rho_f e^V \langle \nabla V \rangle^2 ) \Vert_{L^2(e^{-V/2})}.
\eear
By Lemma \ref{L23}, we estimate
\bear
\nonumber
T_{2,2} &\lesssim& \Vert \eta_1 \Vert_{L^2(e^{V/2})} \Vert \rho_f e^V \langle \nabla V \rangle^2  \Vert_{L^2(e^{-V/2} \langle \nabla V \rangle^{-1}  )} 
\\ \nonumber
&&+ \Vert \eta_2 \Vert_{L^2(e^{V/2})} \Vert \rho_f e^V \langle \nabla V \rangle^2 \Vert_{L^2(e^{-V/2} \langle \nabla V \rangle^{-1} )}
\\ \nonumber
&\lesssim& \Vert f^{\bot} \Vert_\HH \Vert \pi f \Vert_{\HH_1}.
\eear
Using
\bear
\nonumber
j[-\SS f] = j[-\SS f^{\bot}]  &=&- \int v[\Delta_v f^{\bot} + div_v(v f^{\bot})] dv
\\ \nonumber
&=&d\int f^{\bot} v dv  \lesssim \Vert f^{\bot} \Vert_{L^2(e^{v^2/4})},
\eear
and Lemma \ref{L23}, we have
\bear
\nonumber
T_{2,3} &=&(j[-Sf], \nabla (\Delta_V^*)^{-1}(\rho_f e^V \langle \nabla V \rangle^2 ))_{L^2}
\\ \nonumber
&\le &\Vert j[-Sf] \Vert_{L^2( e^{V/2})}  \Vert \nabla (\Delta_V^*)^{-1}(\rho_f e^V \langle \nabla V \rangle^2) \Vert_{L^2( e^{-V/2})}
\\ \nonumber
&\lesssim&\Vert f^{\bot} \Vert_{\HH}  \Vert\rho_f e^V \langle \nabla V \rangle^2 \Vert_{L^2( \langle \nabla V \rangle^{-1} e^{-V/2}) }
\\ \nonumber
&= &\Vert f^{\bot} \Vert_{\HH} \Vert \rho_f  \Vert_{L^2( \langle \nabla V \rangle e^{V/2}  )}
\\ \nonumber
&\lesssim&\Vert f^{\bot} \Vert_{\HH} \Vert \pi f \Vert_{\HH_1}.
\eear
Finally we come to the $T_3$ term. Using 
\beqn
\nonumber
\rho[-Sf] = \int \nabla_v \cdot (\nabla_v f + v f) dv =0,
\eeqn
and
\bear
\nonumber
\rho[-Tf] &= &\rho[v \nabla_x f - \nabla_x V(x) \nabla_v f]
\\ \nonumber
&= &\int v \nabla_x f - \nabla_x V(x) \nabla_v f dv
\\ \nonumber
&= &\nabla_x j[f],
\eear
because $\nabla  (\langle \nabla V \rangle^2) \lesssim \langle \nabla V \rangle^2$ and $ \langle \nabla V \rangle^2 \lesssim \langle \nabla V \rangle$, we get
\bear
\nonumber
T_3 &=&( (\Delta_V)^{-1}\nabla_x j_f, \rho[-\LL f])_{L^2(e^{V/2} \langle \nabla V \rangle)}
\\ \nonumber
&=&((\Delta_V)^{-1} \nabla_x j[f^{\bot}], \rho[-\TT f])_{L^2(e^{V/2}  \langle \nabla V \rangle)}
\\ \nonumber
&=&(j[-f^{\bot}], \nabla (\Delta_V^*)^{-1} (\nabla_x j[f]e^V\langle \nabla V \rangle^2  )_{L^2}
\\ \nonumber
&=&\Vert j[f^{\bot}] \Vert_{L^2(e^{V/2})} \Vert \nabla(\Delta_V^*)^{-1} [\nabla_x(j_f e^V\langle \nabla V \rangle^2)
\\ \nonumber
&&-\nabla V j_f e^V \langle \nabla V \rangle^2 -\nabla (\langle \nabla V \rangle^2) j_f e^V ] \Vert_{L^2(e^{-V/2})},
\eear
using again  Lemma \ref{L23}, we have
\bear
\nonumber
T_3 &\lesssim& \Vert j[f^{\bot}] \Vert_{L^2(e^{V/2})} (\Vert j_f e^V\langle \nabla V \rangle^2 \Vert_{L^2(e^{-V/2} \langle \nabla V \rangle^{-1} )}   
\\ \nonumber
&&+ \Vert j_f e^V \nabla (\langle \nabla V \rangle^2) \Vert_{L^2(\langle \nabla V \rangle^{-1}  e^{-V/2}) })
\\ \nonumber
&\lesssim&\Vert f^{\bot} \Vert_\HH \Vert f \Vert_{\HH_1}.
\eear
Putting all the terms together and choosing $\epsilon > 0$ small enough, we can deduce
\bear
\nonumber
D[f] &\ge& k_p \Vert f^{\bot} \Vert^2_{\HH}+ \epsilon \Vert \pi f \Vert^2_{\HH_1} - \epsilon 2 K\Vert f^{\bot} \Vert_\HH \Vert f\Vert_{\HH_1}- \epsilon 2 K \Vert f^{\bot} \Vert_{\HH} \Vert \pi f \Vert_{\HH_1}
\\ \nonumber
&\ge& k_p \Vert f^{\bot} \Vert^2_{\HH} + \epsilon \Vert \pi f \Vert^2_{\HH_1} - (2\epsilon + 4\epsilon^{1/2}) K \Vert f^{\bot} \Vert^2_{\HH} - \epsilon^{3/2} 4 K \Vert \pi f \Vert^2_{\HH_1}
\\ \nonumber
&\ge& \frac {k_p} 2 (\Vert f^{\bot} \Vert^2_{\HH} + \epsilon \Vert \pi f \Vert^2_{\HH_1}) \ge \frac \epsilon M \Vert f \Vert_{\HH_1},
\eear
for some $M>0$.
\qed

\bigskip
\section{ $L^2$ sub-exponential decay for the kinetic Fokker-Planck equation based on a splitting trick}
\setcounter{equation}{0}
\setcounter{theo}{0}

In this section we establish a first decay estimate on $S_\LL$ which is a particular case in the result of Theorem \ref{T11}.
\begin{theo}\label{T31}
Using the notation and results in Theorem \ref{T21}, we have
\beqn
\nonumber
\Vert S_\LL(t)	f_0\Vert_{L^2(G^{-\frac 1 2})}  \lesssim e^{- C t^{\gamma/(2-\gamma) }} \Vert f_0 \Vert_{L^2 (G^{-(\frac 1 2+\epsilon) } ) } ,
\eeqn
for any $ f_0 \in L^2 (G^{-(\frac 1 2 + \epsilon) } ) \cap \HH_0$,  $   \epsilon > 0$ small enough.
\end{theo}

\begin{rem}
It's worth emphasizing that we deduce immediately part (1) of Theorem \ref{T11} in the case $p=2$ by considering the initial datum $f_0 - \MM(f_0) $ for any $f_0 \in L^2(G^{-\frac {1}  {2} +\epsilon } )$.
\end{rem}

Recall the splitting  $\LL =\AA +\BB$ introduced in (\ref{21}), we first prove some decay estimate on the semigroup $S_\BB$.

\begin{lem} \label{L31} Let us fix $p \in [1,\infty)$. 

(1) For any given smooth weight function $m$, we have
\bear\label{E31}
\int f^{p-1}(\LL f) G^{-(p-1)} m = \frac 1 p \int  f^p G^{-(p-1)} \tilde{m},
\eear
with 
\beqn
\nonumber
\tilde{m} = {\Delta_v m-\nabla_v m \cdot v-\nabla V(x) \cdot \nabla_v m + v \cdot \nabla_x m}  .
\eeqn

(2) Taking $m =e^{\epsilon H^{\delta}} $,  $\epsilon>0$  if $0 < \delta <\frac  \gamma 2$, $\epsilon$ small enough if $\delta =\frac  \gamma 2$, $H = 3v^2 +2 x\cdot v  + x^2 +1$, we have
\bear\label{E32}
\int f^{p-1} (\BB f) G^{-(p-1)} e^{\epsilon H^{\delta}} \le -C \int  f^{p} G^{-(p-1)}e^{\epsilon H^{\delta}} H^{\frac \delta 2 +\gamma -1},
\eear
for some $K$ and $R$ large.

(3) With the same notation as above, there holds

\beqn\label{E33}
\Vert S_\BB(t) \Vert_{L^p(e^{2\epsilon H^\delta}G^{-\frac {p-1} {p}})   \to L^p( e^{\epsilon H^{\delta}} G^{-\frac {p-1} {p}} ) }  \lesssim e^{-at^{\frac {2\delta} {2-\gamma}}},
\eeqn

for some $a>0$. In particular, this implies 

\beqn
\nonumber
\Vert S_\BB(t) \Vert_{L^p(G^{-(\frac {p-1} {p}+ \epsilon) })   \to L^p( G^{-\frac {p-1} {p}} ) }  \lesssim e^{-at^{\frac {\gamma} {2-\gamma}}}.
\eeqn

\end{lem}
\noindent {\sl Proof of Lemma \ref{L31}.}  Step 1. Recall (\ref{A1}), we write
\bear
\nonumber
\int f^{p-1} (\LL f) G^{-(p-1)} m=\int f^{p-1} (\TT f) G^{-(p-1)} m+ \int f^{p-1} (\SS f) G^{-(p-1)} m .
\eear
We first compute the contribution of the term with operator $\TT$
\bear
\nonumber
\int f^{p-1} (\TT f) G^{-(p-1)} m &=& \frac 1 p \int \TT (f^p) G^{-(p-1)} m
\\ \nonumber
&=& -\frac 1 p \int f^p \TT(G^{-(p-1)} m )
\\  \nonumber
&=& \frac 1 p \int f^p G^{-(p-1)} (v \cdot \nabla_x m -\nabla V(x) \cdot \nabla_v m).
\eear
For the the term with operator $\SS$ , we use one integration by parts, and we get
\bear
\nonumber
&&\int f^{p-1} (\SS f) G^{-(p-1)} m
\\ \nonumber
&=& \int f^{p-1} (\Delta_v f +\hbox{div}_v(v f)) G^{-(p-1)} m
\\ \nonumber
&=& -\int \nabla_v((f G^{-1})^{p-1} m ) \cdot  \nabla_v(f G^{-1})G
\\ \nonumber
&=&- \int (p-1) |\nabla_v(f G^{-1})  |^2 (fG^{-1})^{p-2} G m-\frac 1 p  \nabla_v( (f G^{-1})^p )  \cdot (\nabla_v m)  G.
\eear
Performing  another integration by parts on the latter term, we have
\bear
\nonumber
&&\int f^{p-1} (\SS f) G^{-(p-1)} m
\\ \nonumber
&=&\int - (p-1)|\nabla_v(f G^{-1})|^2 (fG^{-1})^{p-2} G m +\frac 1 p \nabla_v \cdot (G \nabla_v m ) (fG^{-1})^{p}
\\ \nonumber
&=&  \int - (p-1)|\nabla_v(f G^{-1} )|^2 (fG^{-1})^{p-2} G m+\frac 1 p (\Delta_v m -v   \cdot \nabla_v  m )f^pG^{-(p-1)}.
\eear
Identity (\ref{E31}) follows by putting together the two identities.\\
Step 2. We particular use $m=e^{\epsilon H^\delta}$ and we easily compute
\beqn
\nonumber
\frac { \nabla_v m} {m}= \delta \epsilon \frac {\nabla_v H} {H^{1-\delta}},\quad \frac {\nabla_x m} {m}=  \delta \epsilon \frac {\nabla_x H} {H^{1-\delta}}  ,
\eeqn
and
\beqn
\nonumber
\frac {\Delta_v m} {m} \le  \delta \epsilon \frac {\Delta_v H} {H^{1-\delta}}+ (\delta \epsilon)^2\frac {|\nabla_v H|^2} {H^{2(1-\delta)}}.
\eeqn
We deduce that $\phi = \frac {\tilde{m}} m$ satisfies 
\bear
\nonumber
&&\frac {\phi H^{1-\delta}} {\epsilon \delta} \le \Delta_v H + \epsilon \delta \frac {|\nabla_v H|^2} {H^{1-\delta}} - v \cdot  \nabla_v H + v \cdot \nabla_x H -\nabla_x V(x) \cdot \nabla_v H.
\eear
From the very definition of $H$, we have
\beqn
\nonumber
\nabla_v H=6 v+ 2 x ,\quad \nabla_x H = 2 v + 2 x,\quad \Delta_v H= 6.
\eeqn
Choosing $\epsilon>0$ arbitrary if $0 < 2\delta < \gamma$, $\epsilon$ small enough if $2\delta = \gamma$ ,we deduce
\bear\label{L2C2}
\nonumber
&&\Delta_v H +2 \epsilon \delta \frac {|\nabla_v H|^2} {H^{1-\delta}}  +v \cdot \nabla_x H - v \cdot \nabla_v H -\nabla_x V(x) \cdot \nabla_v H
\\ \nonumber
&=& 6+ \epsilon \delta \frac{(6 v+2 x)^2} {H^{1-\delta}} + 2 v^2+2 x \cdot  v - 6 v^2- 2 x \cdot v - 6 v \cdot \nabla_x V(x) - 2 x \cdot \nabla_x V(x)
\\ \nonumber
&\le&(2 v^2 +C_1 v + C_2 v^{2\delta} - 6 v^2 ) +(C_3 \epsilon \delta x^{2\delta}- 2 x \cdot \nabla_x V(x)) +C
\\ \nonumber
&\le&-C_4 v ^2 - C_5 x \cdot \nabla_x V(x) + C_6
\\ \nonumber
&\le&- C_7H^{\frac \gamma 2} + K\chi_{R},
\eear
for some constants $C_i, K ,R>0$. As a consequence, we have proved
\beqn
\nonumber
\phi - K\chi_R\le \frac {-C} {H^{1-\delta-\frac \gamma 2}} \le 0 ,
\eeqn
which is nothing but (\ref{E32}).\\
Step 3. In the following, denote $f_t =S_\BB(t)f_0$ the solution to the evolution equation $\partial_t f =\BB f ,f(0)=f_0$. On the one hand, by (\ref{E32}) we have 
\bear
\nonumber
\frac {d} {dt} \int f_t^p G^{-(p-1)} e^{2\epsilon H^{\delta} } = \int f_t^{p-1} (\BB f_t)G^{-(p-1)} e^{2\epsilon H^\delta} \le 0,
\eear
which implies
\beqn
\nonumber
 \int f_t^p G^{-(p-1)} e^{2\epsilon H^{\delta}} \le \int f_0^p G^{-(p-1)} e^{2\epsilon H^{\delta}}:= Y_1, \quad \forall t \ge 0
\eeqn
On the other hand, defining
\beqn
\nonumber
Y :=\int f_t^p G^{-(p-1)} e^{\epsilon H^{\delta} } ,
\eeqn
using again (\ref{E32}),  we have
\bear
\nonumber
\frac {d} {dt} Y&= &p\int f_t^{p-1} \BB f_t G^{-(p-1)}  e^{\epsilon H^{\delta} }  
\\ \nonumber
&\le& -a \int f_t^p G^{-(p-1)} e^{\epsilon H^{\delta} } H^{\delta + \frac \gamma 2 -1}
\\ \nonumber
&\le& -a \int f_t^p G^{-(p-1)}e^{\epsilon H^{\delta} }  \langle x \rangle^{2\delta +\gamma-2}
\\ \nonumber
&\le& -a \int_{B_{|x| \le \rho}} f_t^p G^{-(p-1)} e^{\epsilon H^{\delta} }  \langle x \rangle^{2\delta +\gamma-2},
\eear
for any $\rho > 0$ and for some $a>0$. As $2\delta + \gamma < 2$,  $ 0 \le |x| \le \rho$ implies $\langle x \rangle^{2\delta +\gamma-2} \ge \langle \rho  \rangle^{2\delta +\gamma-2}$, we deduce
\bear
\nonumber
\frac {d} {dt} Y& \le &-a \langle \rho \rangle^{2\delta +\gamma-2}  \int_{B_{|x| \le \rho}} f_t^p G^{-(p-1)}e^{\epsilon H^{\delta} } 
\\ \nonumber
&\le& -a \langle \rho \rangle^{2\delta +\gamma-2} Y + a \langle \rho \rangle^{2\delta +\gamma-2} \int_{B_{|x|  \ge \rho}} f_t^p G^{-(p-1)} e^{\epsilon H^{\delta} } ,
\eear
Using that $e^{\epsilon \langle x \rangle^{2\delta}} \ge e^{\epsilon \langle \rho \rangle^{2\delta}}$ on $ |x| \ge \rho$, we get 
\bear
\nonumber
\frac {d} {dt} Y&\le&  -a \langle \rho \rangle^{2\delta +\gamma-2} Y + a \langle \rho \rangle^{2\delta +\gamma-2} e^{-\epsilon \langle \rho \rangle^{2\delta} } \int_{B_{|x|  \ge \rho}} f_t^p G^{-(p-1)} e^{\epsilon H^{\delta} }  e^{\epsilon \langle x \rangle^{2\delta}}
\\ \nonumber
&\le& -a \langle \rho \rangle^{2\delta +\gamma-2} Y + a \langle \rho \rangle^{2\delta +\gamma-2} e^{-\epsilon \langle \rho \rangle^{2\delta} }  \int f_t^p G^{-(p-1)} e^{\epsilon H^{\delta} } e^{\epsilon \langle x \rangle^{2\delta} }
\\ \nonumber
&\le& -a \langle \rho \rangle^{2\delta +\gamma-2} Y + a \langle \rho \rangle^{2\delta +\gamma-2} e^{-\epsilon  \langle \rho \rangle^{2\delta}}  C Y_1.
\eear 
Thanks to Gr\"{o}nwall's Lemma, we obtain 
\bear
\nonumber
Y(t) &\le& e^{-a \langle \rho \rangle^{2\delta +\gamma-2}t}Y(0) + C e^{- \epsilon \langle \rho \rangle^{2\delta}}Y_1
\\ \nonumber
&\lesssim& (e^{-a \langle \rho \rangle^{2\delta +\gamma-2}t}+ e^{- \epsilon \langle \rho \rangle^{2\delta} })Y_1,
\eear
Choosing finally $\rho$ such that $ a \langle \rho \rangle^{2\delta +\gamma-2} t = \epsilon \langle \rho \rangle^{2\delta}$ , that is $\langle \rho \rangle^{2-\gamma}=C t$, we deduce
\beqn
\nonumber
Y(t) \le C_1 e^{- C_2 t^{\frac {2\delta} {2-\gamma}} }Y_2 ,
\eeqn
for some $C_i>0$, and we deduce the proof of  (\ref{E33}).
\qed 
\\
\sk
\sk
Now we come to prove Theorem \ref{T31}.\\
\noindent {\sl Proof of Theorem \ref{T31}.}
We recall that from (\ref{23}), we have
\beqn
\nonumber
\Vert S_\LL(t) \Vert_{L^2(G^{-1/2} ) \to L^2(G^{-1/2} )} \lesssim  1, \quad \forall t \ge 0
\eeqn
From the very definition of $\AA$ we have
\beqn
\nonumber
\Vert \AA \Vert_{L^2(G^{-1/2} ) \to L^2( e^{2\epsilon H^{\delta}}  G^{-1/2} )} \lesssim 1.
\eeqn
From Lemma \ref{L31} case $p=2$, we have
\beqn
\nonumber
\Vert S_\BB(t) \Vert_{L^2 ( e^{2\epsilon H^{\delta}}  G^{-1/2} ) \to L^2(e^{\epsilon H^{\delta}}  G^{-1/2} )} \lesssim e^{-at^{\frac {2\delta }  { 2-\gamma} }},   \quad \forall t \ge 0.
\eeqn
Gathering the three estimates and using Duhamel's formula 
\beqn
\nonumber
S_\LL = S_\BB + S_\BB \AA * S_\LL
\eeqn
we deduce
\beqn
\nonumber
\Vert S_\LL(t) \Vert_{L^2(e^{2\epsilon H^{\delta}}   G^{-1/2} ) \to L^2(e^{\epsilon H^{\delta}}  G^{-1/2} )} \lesssim 1,   \quad \forall t \ge 0.
\eeqn
In the following, we denote $f_t =S_\LL(t)f_0$ the solution to the evolution equation $\partial_t f =\LL f ,f(0, \cdot)=f_0$. Taking $2\delta =\gamma$, $\epsilon$ small enough, we have in particular
\beqn
\nonumber
\int f_t^2 G^{-1} e^{\epsilon H^{\frac \gamma 2} }  \le C \int f_0^2  G^{-1} e^{2\epsilon H^{\frac \gamma 2} } =:Y_3.
\eeqn
We define
\beqn
\nonumber
Y_2(t):=((f, f)),
\eeqn
with $((, ))$ is defined in Theorem \ref{T21}.
Thanks to the result in (\ref{22}), we have
\bear
\nonumber
\frac {d} {dt} Y_2 &\le& -a \int f_t^2 G^{-1}\langle x \rangle^{2(\gamma-1)}
\\ \nonumber
&\le&-a \int_{B_{|x| \le \rho}} f_t^2 G^{-1} \langle x \rangle^{2(\gamma-1)},
\eear
for any $\rho \ge 0$, using the same argument as Lemma \ref{L31}, we deduce
\bear
\nonumber
Y_2(t) &\le& C e^{-a \langle \rho \rangle^{2(\gamma-1)}t}Y_2(0) + C e^{- \epsilon_2 \langle \rho \rangle^{\gamma}} Y_3
\\ \nonumber
&\lesssim&  (e^{-a \langle \rho \rangle^{2(\gamma-1)}t} + e^{- \epsilon_2 \langle \rho \rangle^{\gamma}})Y_3.
\eear
Choosing $\rho$ such that $ a \langle \rho \rangle^{2(\gamma-1)} t = \epsilon_2 \langle \rho \rangle^{\gamma}$ , that is $\langle \rho \rangle^{2-\gamma}=C t$, we conclude
\beqn
\nonumber
Y_2(t) \le C_1 e^{- C_2 t^{\gamma/(2-\gamma)}} Y_3,
\eeqn
for some constants $C_i>0$. As $H^{\frac \gamma 2} \lesssim C(\frac {v^2} 2 +V(x)) $, we have
\beqn
\nonumber
e^{\epsilon H^{\frac \gamma 2} } \le G^{-C \epsilon},
\eeqn
Taking $\epsilon$ small, the proof of Theorem \ref{T31} is done.
\qed

\bigskip
\section{Regularization property of  $S_\BB$}
\label{sec3}
\setcounter{equation}{0}
\setcounter{theo}{0}

In this section we will denote $\LL^{*}= \LL^{*}_{G^{- 1 /2} }= \SS- \TT$ be the dual operator of $\LL$ on $L^2({G^{-1/2}})$. In other words, $L^*$ is defined by the identity
\beqn
\nonumber
\int (\LL f) g G^{-1}= \int (\LL^{*}g) f G^{-1}.
\eeqn
for any smooth function $f, g$. We also denote $\BB^{*} = \LL^{*} - K\chi_R$.
The aim of this section is to establish the following regularization property. The proof closely follows the proof of similar results in \cite{H, MM, V}
\begin{theo}\label{T41}
For any $0 \le \delta < 1$, there exist  $\eta > 0$ such that
\beqn
\nonumber
\Vert \SS_\BB(t)f \Vert_{L^2(G^{-1/2(1+\delta)})} \lesssim \frac 1 {t^{\frac {5d+1} {2} } } \Vert f \Vert_{L^1(G^{-1/2(1+\delta)})}, \ \ \forall t \in [0,\eta] .
\eeqn
Similarly, for any $0 \le \delta < 1$, there exist  $\eta > 0$ such that
\beqn
\nonumber
\Vert \SS_{\BB^*}(t)f \Vert_{L^2(G^{-1/2(1+\delta)})} \lesssim \frac 1 {t^{\frac {5d+1} {2} } } \Vert f \Vert_{L^1(G^{-1/2(1+\delta}))}, \ \ \forall t \in [0,\eta] .
\eeqn
As a consequence, for any $0 \le \delta < 1$, there exist  $\eta > 0$ such that
\beqn
\nonumber
\Vert \SS_{\BB}(t)f \Vert_{L^\infty(G^{-1/2})} \lesssim \frac 1 {t^{\frac {5d+1} {2} } } \Vert f \Vert_{L^2(G^{-1/2})}, \ \ \forall t \in [0,\eta] .
\eeqn
\end{theo}

We start with some elementary lemmas.
\begin{lem}\label{L41}
For any $0 \le \delta < 1$, we have
\bear\label{E41}
\nonumber
 \quad \int  (f (\LL g) + g (\LL f)) G^{-(1 + \delta )}&=&-2 \int \nabla_v(fG^{-1} )\cdot \nabla_v(gG^{-1 }) G^{1-\delta}
\\ 
&+&\int (\delta d-\delta(1-\delta) v^2)  f g G^{-(1+\delta)}  
\eear
in particular, this implies
\bear\label{E42}
\nonumber
\int f(\LL f) G^{-(1+\delta)} &=& -\int |\nabla_v (fG^{-1})|^2 G^{1-\delta} + \frac {\delta d} {2} \int f^2 G^{-(1+\delta)}
\\ 
&-& \frac {\delta(1-\delta)} {2}  \int v^2 f^2 G^{-(1+\delta)}, 
\eear
similarly, for any $0 \le \delta < 1$, we have
\bear\label{E43}
\nonumber
 \quad  \int f (\LL f) G^{-(1+\delta)}  &=& -\int |\nabla_v f|^2  G^{-(1+\delta)} + \frac { \delta (1 + \delta)} {2} \int v^2 f^2 G^{-(1+\delta)} 
\\ 
&+&\frac {(2+\delta) d }  {2}   \int  f^2G^{-(1+\delta)}.
\eear
All the equalities remain true when $\LL$ is replaced by $\LL^{*}$.
\end{lem}

\sk

\noindent {Proof of Lemma \ref{L41}.} Recall $\TT (G^{-(1+\delta)}) =0$,  we have
\bear
\nonumber
\int f( \TT g) G^{-(1+\delta )} = \int \TT (fG^{-(1 + \delta)}) g =-\int (\TT f) g G^{-(1 + \delta)},
\eear
which implies
\bear
\nonumber
\int f( \TT g) G^{-(1+\delta )} +\int (\TT f) g G^{-(1 + \delta)}=0.
\eear
 for the  term with operator $\SS$ we have
\bear
\nonumber
\int f (\SS g) G^{-(1 +\delta ) } &=& - \int \nabla_v (fG^{-(1 + \delta)})  \cdot (\nabla_v g +v g) 
\\ \nonumber
&=&- \int (\nabla_v f + (1+\delta)v f ) \cdot (\nabla_v g +v g)G^{-(1 + \delta)}
\\ \nonumber
&=& -\int \nabla_v(fG^{-1}) \cdot \nabla_v(gG^{-1}) G^{1-\delta} 
\\ \nonumber
&&- \int( \delta v^2 f  g +\delta f v \cdot  \nabla_v g ) G^{-(1+\delta)} ,
\eear
using integration by parts
\bear
\nonumber
\int \delta f v \cdot \nabla_v g G^{-(1+\delta)} &=& - \int \delta g\nabla_v \cdot ( v f G^{-(1+\delta)}) 
\\ \nonumber
&= &- \int \delta g v \cdot \nabla_v f  G^{-(1+\delta)} 
\\ \nonumber
&& - \int (\delta d  + \delta (1+\delta) v^2 ) f gG^{-(1+\delta)}  ,
\eear
so we deduce
\bear
\nonumber
&&\int (f (\SS g )+ g (\SS f))G^{-(1 +\delta) }
\\ \nonumber
&=& -2 \int \nabla_v(fG^{-1}) \cdot \nabla_v(gG^{-1}) G^{1-\delta}+  \int  ( \delta d  -\delta(1-\delta)v^2   ) f g G^{-(1+\delta)},
\eear
so (\ref{E41}) and (\ref{E42}) are thus proved by combining the two terms above.
Finally, we compute
\bear
\nonumber
&&\int f \SS f G^{-(1 +\delta ) } 
\\ \nonumber
&=& - \int (\nabla_v f + (1+\delta)v f ) \cdot (\nabla_v f +v f)G^{-(1 + \delta)}
\\ \nonumber
&=& -\int |\nabla_v f|^2  G^{-(1+\delta)} - \int (1+\delta) v^2 f^2 G^{-(1+\delta)} - \int (2+ \delta) f v \cdot\nabla_v f G^{-(1+\delta)}
\\ \nonumber
&=& -\int |\nabla_v f|^2  G^{-(1+\delta)} - \int (1+\delta) v^2 f^2 G^{-(1+\delta)} +\frac {2+ \delta} {2} \int  \nabla_v \cdot (v G^{-(1+\delta)}) f^2
\\ \nonumber
&=& -\int |\nabla_v f|^2  G^{-(1+\delta)} + \frac { \delta (1 + \delta)} {2} \int v^2 f^2 G^{-(1+\delta)} +\frac {(2+\delta) d }  {2}   \int f^2 G^{-(1+\delta)},
\eear
so (\ref{E43}) follows by putting together the above equality with 
\beqn
\nonumber
\int f \TT f G^{-(1+\delta)}=0.
\eeqn
Since the term associated with $\TT$ is 0, by $\LL =\SS+\TT, \LL^{*}= \SS - \TT$, we know the same equalities will remain true when $\LL$ is replaced by $\LL^*$.\qed
\sk

\begin{lem}\label{L42}
When $f_t =S_\BB(t) f_0$,  define an energy functional
\bear\label{E44}
\nonumber
\FF(t, f_t ) &:= &A \Vert f_t \Vert_{L^2(G^{-1/2 (1+\delta )})}^2 + at^2 \Vert \nabla_v f_t \Vert_{L^2(G^{-1/2(1+\delta) ) })}^2 
\\ 
&+& 2 c t^4( \nabla_v f_t, \nabla_x f_t )_{L^2(G^{-1/2  (1+\delta) })}^2 + bt^6 \Vert \nabla_x f_t \Vert_{L^2(G^{-1/2 (1+\delta) )})}^2,
\eear
when $f_t =S_{\BB^*}(t) f_0$,  define another energy functional
\bear\label{E45}
\nonumber
\FF^*(t, f_t ) &:= &A \Vert f_t \Vert_{L^2(G^{-1/2 (1+\delta )})}^2 + at^2 \Vert \nabla_v f_t \Vert_{L^2(G^{-1/2(1+\delta) ) })}^2 
\\
&-& 2 c t^4 ( \nabla_v f_t, \nabla_x f_t )_{L^2(G^{-1/2  (1+\delta) })}^2 + bt^6 \Vert \nabla_x f_t \Vert_{L^2(G^{-1/2 (1+\delta) )})}^2,
\eear
with $a, b, c >0, c \le \sqrt{ab} $ and $A$ large enough. Then for both cases we have
\bear
\nonumber
\frac {d } {dt}F(t, f_t )\le -L(\Vert \nabla_v f_t \Vert_{L^2(G^{-1/2(1+\delta )})}^2 + t^4 \Vert \nabla_x f _t \Vert_{L^2 ( G^{-1/2  (1+\delta)}) }^2 ) +\Vert f_t \Vert_{L^2(G^{-1/2(1+\delta )})}^2 ,
\eear
for all $t \in [0, \eta]$, for some $L>0, C>0$ and $F= \FF$ or $\FF^*$.
\end{lem}

\noindent {Proof of Lemma \ref{L42}.} We only prove the case $F=\FF$, the proof for $F=\FF^*$ is the same. We split the computation into several parts and then put them together. First using (\ref{E42}) and (\ref{E43}) we have
\bear
\nonumber
&&\frac {d} {dt} \int f^2 G^{-(1+\delta)}  
\\ \nonumber
&= &\int f (\LL -K\chi_R) f G^{-(1+\delta)}  
\\ \nonumber
&=&\frac {1-\delta} 2 \int f \LL f  G^{-(1+\delta)} +\frac {1+\delta} 2\int f \LL f  G^{-(1+\delta)} -\int K\chi_R f^2 G^{-(1+\delta)}
\\ \nonumber
&\le& -\frac {1-\delta} 2\int |\nabla_v f|^2  G^{-(1+\delta)} -\frac {1+\delta} 2\int |\nabla_v (fG^{-1})|^2 G^{1-\delta} + C \int f^2 G^{-(1+\delta)}
\\ \nonumber
&\le& -\frac {1-\delta} 2\int |\nabla_v f|^2  G^{-(1+\delta)} + C \int f^2 G^{-(1+\delta)}.
\eear
By
\beqn\label{E46}
\partial_{x_i}\LL f =  \LL \partial_{x_i} f  + \sum\limits_{j=1}^d \partial^2_{x_i x_j} V\partial_{v_j} f,
\eeqn  
and (\ref{E42}) we have
\bear
\nonumber
&&\frac {d} {dt} \int (\partial_{x_i} f)^2 G^{-(1+\delta)} 
\\ \nonumber
&=&\int \partial_{x_i} f \partial_{x_i} (\LL- K\chi_R) f G^{-(1+\delta)}
\\ \nonumber
&=& -\int |\nabla_v (\partial_{x_i} f G^{-1})|^2 G^{1-\delta} +  \frac {\delta d} {2} \int (\partial_{x_i} f)^2 G^{-(1+\delta)} 
\\ \nonumber
&&-\frac {\delta(1-\delta) } {2}  \int v^2 (\partial_{x_i} f)^2 G^{-(1+\delta)}+ \int \partial_{x_i} f\sum\limits_{j=1}^d \partial^2_{x_i x_j} V\partial_{v_j} f G^{-(1+\delta)}
\\ \nonumber
&&-\int  K\chi_R |\partial_{x_i} f|^2 G^{-(1+\delta)} -   \int K \partial_{x_i} f \partial_{x_i}\chi_R  f G^{-(1+\delta)}.
\eear
Using Cauchy-Schwarz inequality and summing up by i, we get
\bear
\nonumber
&&\frac {d} {dt} \int |\nabla_x f|^2 G^{-(1+\delta)}  
\\ \nonumber
&\le& -\sum\limits_{i=1}^d \int |\nabla_v (\partial_{x_i} f G^{-1})|^2 G^{1-\delta} - \frac {\delta(1-\delta)} {2}  \int v^2 (\nabla_x f)^2 G^{-(1+\delta)}
\\ \nonumber
&&+C\int |\nabla_v f |^2 G^{-(1+\delta)} +C \int |\nabla_{x} f|^2 G^{-(1+\delta)}+ C \int | f|^2 G^{-(1+\delta)} ,
\eear
for some $C > 0$. Similarly using
\beqn\label{E47}
\partial_{v_i}\LL f=\LL\partial_{v_i} f -\partial_{x_i}f +\partial_{v_i} f ,
\eeqn
and (\ref{E42}), we have
\bear
\nonumber
&&\frac {d} {dt} \int (\partial_{v_i} f)^2 G^{-(1+\delta)}  
\\ \nonumber
&=&\int \partial_{v_i} f \partial_{v_i} (\LL-K\chi_R) f G^{-(1+\delta)}
\\ \nonumber
&=& -\int |\nabla_v (\partial_{v_i} f G^{-1})|^2 G^{1-\delta} + \frac {\delta d} {2} \int  (\partial_{v_i} f)^2 G^{-(1+\delta)} 
\\ \nonumber
&&-\frac {\delta(1-\delta)} {2}  \int v^2  (\partial_{v_i} f)^2 G^{-(1+\delta)}- \int \partial_{x_i} f \partial_{v_i} f G^{-(1+\delta)}
\\ \nonumber
&&+ \int |\partial_{v_i} f|^2 G^{-(1+\delta)}- \int  K\chi_R |\partial_{v_i} f|^2 G^{-(1+\delta)}- \int K \partial_{v_i} f \partial_{v_i}\chi_R  f G^{-(1+\delta)}.
\eear
Using Cauchy-Schwarz inequality and summing up by i we get
\bear
\nonumber
&&\frac {d} {dt} \int |\nabla_v f|^2 G^{-(1+\delta)}  
\\ \nonumber
&\le&- \sum\limits_{i=1}^d \int |\nabla_v (\partial_{v_i} f G^{-1})|^2 G^{1-\delta} + C\int |\nabla_x f||\nabla_v f|G^{-(1+\delta)}
\\ \nonumber
&& +C\int |\nabla_{v} f|^2 G^{-(1+\delta)} +   C\int  |f|^2 G^{-(1+\delta)} - \frac {\delta(1-\delta)} {2}  \int v^2  (\nabla_v  f)^2 G^{-(1+\delta)}.
\eear
For the crossing term, we split it also into two parts 
\bear
\nonumber
&&\frac d {dt} \int 2\partial_{v_i} f \partial_{x_i} f G^{-(1+\delta)}  
\\ \nonumber
&=& (\int \partial_{v_i} f \partial_{x_i} \LL f  G^{-(1+\delta)}  +\int \partial_{v_i} \LL f \partial_{x_i}  f  G^{-(1+\delta)})
\\ \nonumber
&&- (\int \partial_{v_i} f \partial_{x_i}(K\chi_R f) G^{-(1+\delta)} +\int \partial_{x_i} (K\chi_R f) \partial_{v_i } f G^{-(1+\delta)})
\\ \nonumber
&:=& W_1 +W_2.
\eear
Using (\ref{E46}) and (\ref{E47}) we have
\bear
\nonumber
W_1&=& \int \partial_{v_i} f \LL (\partial_{x_i} f) G^{-(1+\delta)} +\int \LL(\partial_{v_i}  f) \partial_{x_i}  f G^{-(1+\delta)}
\\ \nonumber
&&+ \int\partial_{v_i} f \sum\limits_{j=1}^d  \partial_{x_i x_j} V(x) \partial_{v_j} f G^{-(1+\delta)} -\int |\partial_{x_i}f|^2 G^{-(1+\delta)} 
\\ \nonumber
&&+\int  \partial_{x_i} f \partial_{v_i} f G^{-(1+\delta)}.
\eear
By (\ref{E41}), we deduce
\bear
\nonumber
W_1&=&- \int 2 \nabla_v(\partial_{v_i}fG^{-1}) \cdot \nabla_v(\partial_{x_i} fG^{-1}) G^{1-\delta}+  \delta d \int \partial_{v_i} f \partial_{x_i} f G^{-(1+\delta)} 
\\ \nonumber
&&-\delta(1-\delta)  \int v^2 \partial_{v_i} f \partial_{x_i} f G^{-(1+\delta)}+ \int\partial_{v_i} f \sum\limits_{j=1}^d  \partial_{x_i x_j} V(x) \partial_{v_j} f G^{-(1+\delta)}
\\ \nonumber
&&-\int |\partial_{x_i}f|^2 G^{-(1+\delta)} +\int  \partial_{x_i} f \partial_{v_i} f G^{-(1+\delta)}.
\eear
For the $W_2$ term we have
\bear
\nonumber
W_2&=&- \int \partial_{v_i} f \partial_{x_i}(K\chi_R f) G^{-(1+\delta)} -\int \partial_{x_i} (K\chi_R f) \partial_{v_i } f G^{-(1+\delta)}
\\ \nonumber
 &=& - \int 2 K\chi_R \partial_{x_i} f \partial_{v_i} f  G^{-(1+\delta)}  + \int K f  (\partial_{v_i} \chi_R \partial_{x_i} f+ \partial_{v_i} f \partial_{x_i} \chi_R)G^{-(1+\delta)}
\\ \nonumber
&\le&  C\int |\partial_{x_i} f |  | \partial_{v_i} f | G^{-(1+\delta)} + C \int |\partial_{v_i} f || f | G^{-(1+\delta)}+ C\int |f | |\partial_{x_i} f | G^{-(1+\delta)},
\eear
Combining the two parts, using Cauchy-Schwarz inequality, and summing up by i we get
\bear
\nonumber
&&\frac {d} {dt} \int 2  \nabla_x f \cdot \nabla_v f G^{-(1+\delta)}
\\ \nonumber
&\le&  - \sum\limits_{i=1}^d  \int 2 \nabla_v(\partial_{v_i}fG^{-1}) \cdot \nabla_v(\partial_{x_i} fG^{-1}) G^{1-\delta} -\frac 1 2 \int |\nabla_x f|^2 G^{-(1+\delta)}
\\ \nonumber
&& + C \int |\nabla_v f|^2 G^{-(1+\delta)} +  C\int |f |^2 G^{-(1+\delta)} -\delta(1-\delta)  \int v^2 \nabla_v f \cdot \nabla_x f G^{-(1+\delta)}.
\eear
For the very definition of $\FF$ in (\ref{E44}), we easily compute
\bear
\nonumber
&&\frac d {dt} \FF(t, f_t) 
\\ \nonumber
&=& A \frac d {dt} \Vert f_t \Vert_{L^2(G^{-1/2(1+\delta)})}^2 + at^2   \frac d {dt} \Vert \nabla_v f_t \Vert_{L^2(G^{-1/2(1+\delta)})}^2 
\\ \nonumber
&&+2 c t^4 \frac d {dt} \langle \nabla_v f_t, \nabla_x f_t \rangle_{L^2(G^{-1/2(1+\delta)})}^2 + bt^6 \frac d {dt} \Vert \nabla_x f_t \Vert_{L^2(G^{-1/2(1+\delta)})}^2
\\ \nonumber
&&+ 2 a t \Vert \nabla_v f_t \Vert_{L^2(G^{-1/2  (1+\delta) })}^2 +  8c t^3 \langle \nabla_v f_t, \nabla_x f_t \rangle_{L^2(G^{-1/2 (1+\delta) })}^2 
\\ \nonumber
&&+ 6 b t^5 \Vert \nabla_x f_t \Vert_{L^2(G^{-1/2(1+\delta)  })}^2.
\eear
Gathering all the inequalities above together, we have
\bear
\nonumber
&&\frac d {dt} \FF(t, f_t) 
\\ \nonumber
&\le&  (2at - \frac {A(1-\delta)}  2+ Cat^2 + 2Ct^4c +C b t^6) \int | \nabla_v f_t |^2 G^{-(1+\delta)} 
\\ \nonumber
&&+ (6bt^5-  \frac c 2 t^4+C b t^6)\int |\nabla_x f_t|^2 G^{-(1+\delta)} + (8ct^3 +Cat^2 ) \int |\nabla_{v} f_t | |\nabla_{x} f_t | G^{-(1+\delta)}  
\\ \nonumber
&& -( a t^2 \sum\limits_{i=1}^{d} \int |\nabla_v (\partial_{v_i} f_t G^{-1})|^2 G^{1-\delta}  +b t^6 \sum\limits_{i=1}^{d} \int |\nabla_v (\partial_{x_i} f_t G^{-1})|^2 G^{1-\delta}
\\ \nonumber
&& +  2 ct^4 \sum\limits_{i=1}^{d}   \int  \nabla_v(\partial_{v_i} f_t G^{-1}) \cdot \nabla_v(\partial_{x_i} f_t G^{-1}) G^{1-\delta} ) - \frac {\delta(1-\delta)} {2} (at^2 \int v^2  (\nabla_v  f)^2 G^{-(1+\delta)}
\\ \nonumber
&&  + bt^6\int v^2  (\nabla_x  f)^2 G^{-(1+\delta)} + 2c t^4\int v^2 \nabla_v f \cdot \nabla_x f G^{-(1+\delta)} )  +C\int f_t^2 G^{-(1+\delta)},
\eear
for some $C>0$. We observe that
\bear
\nonumber
&&|2c t^4\int v^2 \nabla_v f \cdot \nabla_x f G^{-(1+\delta)} |
\\ \nonumber
&\le& at^2 \int v^2  (\nabla_v  f)^2 G^{-(1+\delta)} + bt^6\int v^2  (\nabla_x  f)^2 G^{-(1+\delta)} ,
\eear
and
\bear
\nonumber
&&| 2 ct^4 \sum\limits_{i=1}^{d}   \int 2 \nabla_v(\partial_{v_i}f_tG^{-1}) \cdot \nabla_v(\partial_{x_i} f_tG^{-1}) G^{1-\delta} |
\\ \nonumber 
&\le &   a t^2 \sum\limits_{i=1}^{d} \int |\nabla_v (\partial_{v_i} f_t G^{-1})|^2 G^{1-\delta}  +b t^6 \sum\limits_{i=1}^{d} \int |\nabla_v (\partial_{x_i} f_t G^{-1})|^2 G^{1-\delta}.
\eear
by our choice on $a, b, c$. So by taking $A $ large and $0 < \eta$ small ($t \in [0, \eta]$), as a consequence
\bear
\nonumber
\frac d {dt} \FF(t, f_t)\le -L(\Vert \nabla_v f_t \Vert_{L^2(G^{-1/2(1+\delta )})}^2 + t^4 \Vert \nabla_x f _t \Vert_{L^2 ( G^{-1/2  (1+\delta)}) }^2) +C\Vert f_t \Vert_{L^2(G^{-1/2(1+\delta )})}^2,
\eear
for some $L, C  > 0$, and that ends the proof.\qed   \sk
\begin{rem}
For the case $F = \FF^*$, the only difference in the proof is to change (\ref{E46}) and (\ref{E47}) into
\beqn
\nonumber
\partial_{x_i}\LL^* f = \LL^* \partial_{x_i} f - \partial_{x_i}(\nabla_x V(x) \cdot \nabla_v f)=  \LL^* \partial_{x_i} f  - \sum\limits_{j=1}^d \partial^2_{x_i x_j} V\partial_{v_j} f,
\eeqn
and
\beqn
\nonumber
\partial_{v_i}\LL^* f=\LL^* \partial_{v_i} f +\partial_{x_i}f +\partial_{v_i} f.
\eeqn
\end{rem}
\sk

The following proof of this section is true for both cases.
\begin{lem}\label{L53} For any $0 \le \delta <1$, we have
\bear
\nonumber
\int |\nabla_{x, v} (f G^{-1/2(1+\delta)}) |^2 \le \int |\nabla_{x, v} f |^2 G^{-(1+\delta)} +C \int f^2 G^{-(1+\delta)},
\eear
\end{lem}
\noindent {Prove of Lemma   \ref {L53}.}
For any weight function $m$ we have
\bear
\nonumber
\int |\nabla_x (f m) |^2 &=& \int|\nabla_x f m +\nabla_x m f |^2 
\\ \nonumber
&=&\int |\nabla_x f |^2 m^2  +\int |\nabla_x m |^2 f^2 +\int 2 f\nabla_x f  m \nabla_x m
\\ \nonumber
&=& \int |\nabla_x f |^2 m^2 +\int (|\nabla_x m |^2-\frac 1 2 \Delta_x (m^2)) f^2,
\eear
taking $m = G^{-1/2(1+\delta)}$ we have
\bear
\nonumber
&&\int |\nabla_x (f G^{-1/2(1+\delta)}) |^2 
\\ \nonumber
&=&\int |\nabla_x f |^2 G^{-(1+\delta)} +\int -(\frac {(1+\delta)^2}  4|\nabla_x V(x)|^2 +\frac {1+\delta} 2 \Delta_x V(x)) f^2 G^{-(1+\delta)}
\\ \nonumber
&\le& \int |\nabla_x f |^2 G^{-(1+\delta)} +C \int f^2 G^{-(1+\delta)}.
\eear
Similarly, we have
\bear
\nonumber
&&\int |\nabla_v (f G^{-1/2(1+\delta)}) |^2 
\\ \nonumber
&=&\int |\nabla_v f |^2 G^{-(1+\delta)} +\int -(\frac {(1+\delta)^2} 4v^2 +\frac {1+\delta} 2 d) f^2 G^{-(1+\delta)}
\\ \nonumber
&\le& \int |\nabla_v f |^2 G^{-(1+\delta)}.
\eear
Putting together the two inequalities we obtain the result.
\qed

\begin{lem}\label{L54}
Nash's inequality: for any $f \in L^1(\R^d) \cap H^{1}(\R^d) $,there exist a constant $C_d$ such that:
\beqn
\nonumber
\Vert f \Vert_{L^2}^{1+\frac 2 d} \le C_d \Vert f \Vert_{L^1}^{ 2/ d}\Vert \nabla_v f \Vert_{L^2},
\eeqn
 \end{lem} 
\noindent For the proof of Nash's inequality, we refer to \cite{LL}, Section 8.13 for instance.
\qed

\begin{lem}\label{L55} For any $0 \le \delta <1$ we have
\bear\label{E48}
\frac d {dt} \int  | f | G^{-1/2(1+\delta)} \le d  \int |f | G^{-1/2(1+\delta)},
\eear
which implies
\bear
\nonumber
\int  | f_t | G^{-1/2(1+\delta)} \le Ce^{ d  t }\int |f_0 | G^{-1/2(1+\delta)}.
\eear
In particular we have
\bear\label{E49}
\int  | f_t | G^{-1/2(1+\delta)} \le C\int |f_0 | G^{-1/2(1+\delta)} , \ \ \ \forall t \in [0,\eta],
\eear
for some constant $C>0$.
\end{lem}

\noindent {Proof of  Lemma \ref{L55}.} By Lemma \ref{L51} in the next section, letting $p=1$, we have
\bear
\nonumber
&&\frac d {dt} \int  | f | G^{-1/2(1+\delta)} 
\\ \nonumber
&=& \int |f | (\Delta_v G^{-1/2(1+\delta)} -v \cdot \nabla_v G^{-1/2(1+\delta)}
\\ \nonumber
&&+ v \cdot \nabla_x G^{-1/2(1+\delta)} - \nabla V(x) \cdot \nabla_v G^{-1/2(1+\delta)} -K\chi_R G^{-1/2(1+\delta)})
\\ \nonumber
&\le& \int |f |(\frac {1+\delta}  2 d  -\frac {(1+\delta) (1-\delta)} 4 v^2) G^{-1/2(1+\delta)} \le \int |f | d G^{-1/2(1+\delta)}.
\eear
so (\ref{E48}) is proved. As $\TT G^{-1/2(1+\delta)} =0$, the result is still true when $F=\FF^*$.
\qed

\sk

Now we come to the proof of Theorem \ref{T41}.

\sk

\noindent {Proof of Theorem \ref{T41}.} We define
\beqn
\nonumber
\GG(t, f_t)=B \Vert  f_t \Vert_{L^1(G^{-1/2(1+\delta)})}^2 + t^Z \FF (t, f_t),
\eeqn
with $B, Z > 0 $ to be fixed and $\FF$ defined in Lemma \ref{L42}. We choose $t \in [0, \eta]$ , $\eta$ small such that $(a+b+c) Z\eta^{Z+1} \le \frac 1 2 L \eta^Z$ ($a, b, c, L$ are also defined Lemma \ref{L42}), by  (\ref{E48}) and Lemma \ref{L42} we have
\bear
\nonumber
\frac  d {dt} \GG(t, f_t) &\le& dB \Vert f_t \Vert_{L^1(G^{-1/2(1+\delta)})}^2 +Zt^{Z-1} \FF(t, f_t) 
\\ \nonumber
&&- L t^Z(\Vert \nabla_v f_t \Vert_{L^2(G^{-1/2(1+\delta)})}^2 + t^4 \Vert \nabla_x f _t \Vert_{L^2 ( G^{-1/2(1+\delta)} ) }^2 )
\\ \nonumber
&&+Ct^Z \Vert  f_t \Vert_{L^2(G^{-1/2(1+\delta )})}^2 
\\ \nonumber
&\le& dB \Vert f_t \Vert_{L^1(G^{-1/2(1+\delta)} ) }^2  +Ct^{Z-1}\Vert  f_t \Vert_{L^2(G^{-1/2(1+\delta )})}^2 
\\ \nonumber
&&- \frac L 2 t^Z(\Vert \nabla_v f_t \Vert_{L^2( G^{-1/2(1+\delta)} )}^2 + t^4 \Vert \nabla_x f _t \Vert_{L^2 ( G^{-1/2(1+\delta)} ) }^2 ).
\eear
Nash's inequality and Lemma \ref{L42} implies
\bear 
\nonumber
\int f_t^2 G^{-(1+\delta)} &\le& (\int |f_t| G^{-1/2(1+\delta)})^{\frac 4 {d+2}} (\int |\nabla_{x, v}(f_t G^{-1/2(1+\delta)}) |^2 )^{\frac d {d+2}}
\\ \nonumber
&\le&  (\int |f_t| G^{-1/2(1+\delta)} )^{\frac 4 {d+2}} (\int |\nabla_{x, v}f_t  |^2 G^{-(1+\delta)} +C\int f_t^2 G^{-(1+\delta)})^{\frac d {d+2}}.
\eear
Using Young's inequality, we have
\beqn
\nonumber
\Vert f_t \Vert_{L^2( G^{-1/2(1+\delta)} )}^2 \le C_{\epsilon} t^{-5d} \Vert f \Vert_{L^1(G^{-1/2(1+\delta)} )}^2 + \epsilon t^5 (\Vert \nabla_{x, v} f_t \Vert_{L^2( G^{-1/2(1+\delta)} )}^2 + C\Vert  f_t \Vert_{L^2(G^{-1/2(1+\delta)} )}^2).
\eeqn
Taking $\epsilon$ small such that $C \epsilon \eta^5 \le  \frac 1 2$, we deduce
\beqn
\nonumber
\Vert f_t \Vert_{L^2(G^{-1/2(1+\delta)})}^2 \le 2C_{\epsilon} t^{-5d} \Vert f \Vert_{L^1(G^{-1/2(1+\delta)})}^2 + 2\epsilon t^5 \Vert \nabla_{x, v} f_t \Vert_{L^2(G^{-1/2(1+\delta)})}^2 .
\eeqn
Taking $\epsilon$ small we have
\beqn
\nonumber
\frac d {dt} \GG(t, f_t) \le dB \Vert f_t \Vert_{L^1(G^{-1/2(1+\delta)})}^2+ C_1t^{Z-1-5d} \Vert f_t \Vert_{L^1(G^{-1/2(1+\delta)})}^2,
\eeqn
for some $C_1>0$. Choosing $Z=1+5d$, and using (\ref{E49}), we deduce
\beqn
\nonumber
\forall t \in [0, \eta], \ \ \ \GG(t, f_t) \le \GG(0, f_0) + C_2\Vert f_0 \Vert_{L^1(G^{-1/2(1+\delta)} )}^2 \le C_3 \Vert f_0 \Vert_{L^1(G^{-1/2(1+\delta)})}^2,
\eeqn
which ends the proof.
\qed

\bigskip
\section{ $S_\BB$ decay in larger spaces}
\label{sec:L1}
\setcounter{equation}{0}
\setcounter{theo}{0}

The aim of this section is to prove the following decay estimate for the semigroup $S_\BB$ which will be useful in the last section where we will prove Theorem \ref{T11} in full generally. 
\begin{theo}\label{T51}
Let $H=1 + x^2+ 2 v \cdot x  +3v^2$, for any $\theta \in(0, 1)$ and
for any $l > 0$, we have 
\beqn
\nonumber
\Vert \SS_\BB(t)  \Vert_{L^1(H^l) \to L^1(H^{l \theta})}\lesssim (1+t)^{-a},
\eeqn
where
\beqn
\nonumber
a = \frac {l(1- \theta)} {1- \frac \gamma 2 }.
\eeqn

\end{theo}

We start with an elementary identity.
\begin{lem}\label{L51}
For the kinetic Fokker Planck operator $\LL$ , let m  be a  weight function, for any $p \in [1, \infty]$ we have 
\beqn
\\ \nonumber
 \int (\LL f)f^{p-1}m^p=-(p-1)\int |\nabla_v(mf)|^2 (mf)^{p-2} + \int f^pm^p \phi ,
 \eeqn
with
\beqn
\nonumber
\phi =  \frac {2} {p^{'}} \frac {|\nabla_v m|^2} {m^2} + (\frac 2 p -1)\frac {\Delta_v m} {m} + \frac {d} {p^{'}}- v \cdot \frac {\nabla_v m} {m} -\frac {\TT m } {m}.
\eeqn
In particular when $p=1$, we have
\beqn
\nonumber
\phi = \frac {\Delta_v m} {m}- v \cdot \frac {\nabla_v m} {m} -\frac {\TT m } {m}.
\eeqn
\end{lem}
\noindent {\sl Proof of Lemma \ref{L51}.} We split the integral as
\beqn
\nonumber
\int (\LL f)f^{p-1}m^p=\int f^{p-1}\SS f m^p  +\int f^{p-1}\TT f m^p.
\eeqn
First compute the contribution of the term with operator $\TT$
\beqn
\nonumber
\int f^{p-1}\TT f m^p =\frac 1 p  \int  \TT (f^p)m^p = -\int f^p m^p \frac {\TT m} {m} .
\eeqn
\smallskip
Concerning the term with operator $\SS$, we split it also into two parts
\beqn
\nonumber
\int (\SS f)f^{p-1}m^p = \int f^{p-1} m^p (\Delta_v f +div_v(v f)) := C_1+C_2.
\eeqn
We first compute the $C_2$ term, we have
\bear
\nonumber
C_2&=&\int f^{p-1} m^p(d f+ v \cdot \nabla_v f)
\\ \nonumber
&=&\int d f^pm^p - \frac 1 p \int f^p \hbox{div}_v(v m^p)
\\ \nonumber
&=&\int f^p[(1-\frac 1 p )d -v\cdot \frac {\nabla_v m}{m}]m^p.
\eear
We turn to the $C_1$ term, we have
\bear
\nonumber
C_1&=& \int f^{p-1} m^p \Delta_v f =- \int \nabla_v (f^{p-1}m^p) \cdot \nabla_v f
\\ \nonumber
&=&\int - (p-1)|\nabla_v f|^2 f^{p-2} m^{p} - \frac 1 p \int \nabla_v f^p  \cdot \nabla_v m^p.
\eear
Using $\nabla_v(m f) =m \nabla_v f + f \nabla_v m$, we deduce
\bear
\nonumber
C_1&=&-(p-1) \int |\nabla_v(mf)|^2f^{p}m^{p-2}+(p-1)\int |\nabla_v m|^2f^p m^{p-2}
\\ \nonumber
&&+ \frac{2(p-1)} {p^2} \int \nabla_v (f^p) \cdot  \nabla_v (m^{p}) - \frac 1 p \int \nabla_v (f^p) \cdot \nabla_v (m^p)
\\ \nonumber
&=&-(p-1) \int |\nabla_v(m f)|^2f^{p-2} m^p+(p-1)\int |\nabla_v m|^2f^p m^{p-2}
\\ \nonumber
&& + \frac {p-2} {p^2} \int f^p \Delta_v m^p.
\eear
Using that $\Delta_v m^p = p \Delta_v m\  m^{p-1} + p(p-1)| \nabla_v m|^2m^{p-2}  $, we obtain
\bear
\nonumber
C_1=-(p-1) \int |\nabla_v(mf)|^2f^{p-2}m^{p-2} + \int f^p m^p[ (\frac 2 p -1) \frac { \Delta_v m} {m} + 2(1-\frac 1 p)\frac {|\nabla_v m|^2} {m^2}].
\eear
We conclude by combining the above equalities.
\qed
\smallskip

\smallskip
\smallskip
\noindent {\sl Proof of Theorem \ref{T51}.} From Lemma \ref{L51}, we have
\bear\label{E51}
 \int (\BB f)f^{p-1}m^p &&= \int (\LL -M \chi_R) f^{p-1} m^p
 \\ \nonumber
&&= - (p-1)\int |\nabla_v(mf)|^2 (mf)^{p-2} + \int f^pm^p \phi ,
 \eear
with
\beqn
\nonumber
\phi =  [\frac {2} {p^{'}} \frac {|\nabla_v m|^2} {m^2} + (\frac 2 p -1)\frac {\Delta_v m} {m} + \frac {d} {p^{'}}- v \cdot \frac {\nabla_v m} {m} -\frac {\TT m } {m} - M\chi_R].
\eeqn
When $p=1$, we have
\bear
\nonumber
\phi =  \frac {\Delta_v m} {m} - v \cdot \frac {\nabla_v m} {m} -\frac {\TT m } {m} - M\chi_R.
\eear
Let $m=H^k$. We have
\bear
\nonumber
\frac {\nabla_v m} {m}= k \frac {\nabla_v H} {H}, \quad \frac {\nabla_x m} {m}=  k \frac {\nabla_x H} {H} ,
\eear
and 
\bear
\nonumber
\frac {\Delta_v m} {m}=  \frac {k\Delta_v H} {H} + \frac {k (k-1)|\nabla_v H|^2} {H^2}.
\eear
Summing up, we have for $\phi$
\bear
\nonumber
\frac {\phi H } {k} =  \Delta_v H + (k-1)\frac {|\nabla_v H|^2} {H} - v \cdot  \nabla_v H + v \cdot \nabla_x H -\nabla_x V(x) \cdot \nabla_v H -M \chi_R,
\eear
From the very definition of $H$, we have 
\beqn 
\nonumber
\nabla_v H= 6 v+ 2 x , \quad \nabla_x H = 2 v + 2 x, \quad \Delta_v H= 6.
\eeqn
We then compute 
\bear\label{L2C2}
\nonumber
&&\Delta_v H + (k-1)\frac {|\nabla_v H|^2} {H}  +v \cdot \nabla_x H - v \cdot \nabla_v H -\nabla_x V(x)\cdot \nabla_v H
\\ \nonumber
&=& 6+ (k-1)\frac{( 6v+2 x)^2} {H} + 2 v^2+2 x \cdot v -  6v^2
\\ \nonumber
&&- 2 x \cdot v -  6 v \cdot \nabla_x V(x) - 2 x \cdot \nabla_x V(x)
\\ \nonumber
&\le&(2 v^2 +C v - 6v^2) - 2 x \cdot \nabla_x V(x) + C
\\ \nonumber
&\le& - C_1 v ^2 - C_2 x \cdot \nabla_x V(x) + C_3
\\ \nonumber
&\le& -C_4 H^{\frac \gamma 2}  +K_1\chi_{R_1} ,
\eear
for some $C_i >0$. Taking $K$ and $R$ large enough, we have $ \phi \le -C H^{\frac \gamma 2 - 1}$, using  this inequality in equation (\ref{E51}), we deduce
\bear\label{E52}
\frac d {dt} Y_4(t) :=\frac d {dt} \int |f_\BB(t)| H^k &=& \int sign (f_\BB(t)) \BB f_\BB(t) H^k
\\ \nonumber
&\le& -C \int |f_B(t)| H^{k -1+ \frac \gamma 2} , 
\eear
for any $k >1$. In particular for any $l \ge 1$, we can find $K$ and $R$ large enough such that
\beqn
\nonumber
\frac d {dt} \int |f_\BB(t)| H^l \le  0,
\eeqn
which readily implies
\beqn
\nonumber
\int |f_\BB(t) |H^l \le \int |f_0| H^l := Y_5.
\eeqn
Denoting 
\beqn
\nonumber
\alpha=\frac {l-k} {l-k+1-\frac \gamma 2} \in (0,1),
\eeqn
the H\"older's inequality
\beqn
\nonumber
\int |f_B(t)| H^{k } \le (\int |f_B(t)| H^{k -1+\frac \gamma 2} )^{\alpha} (\int |f_B(t)| H^{l}  )^{1-\alpha},
\eeqn
implies
\beqn
\nonumber
(\int |f_B(t)| H^{k} )^{\frac 1 \alpha }(\int |f_B(t)| H^{l}  )^{\frac {\alpha-1} {\alpha}} \le \int |f_B(t)| H^{k -1+\frac \gamma 2}  ,
\eeqn
From this inequality and (\ref{E52}), we get
\beqn
\nonumber
\frac d {dt} Y_4(t) \le -C (Y_4(t))^{\frac 1 \alpha} Y_5^{\frac {\alpha-1} {\alpha}}.
\eeqn
Using $Y_4(0) \le Y_5$, after an integration, we deduce
\beqn
\nonumber
Y_4(t) \le C_\alpha \frac 1 {(1+t)^{\frac {\alpha} {1-\alpha}} }Y_5,
\eeqn
which is nothing but the polynomial decay on $S_\BB$
$$\Vert \SS_B(t)  \Vert_{L^p(H^l) \to L^p(H^k)}\lesssim (1+t)^{-a},$$
with
$$ a = \frac {l-k} {1- \frac \gamma 2 }, \quad \forall 0  < k < l, \quad 1\le l.$$
We conclude Theorem \ref{T51} by writing $k = l\theta$, $0<\theta<1$.

\bigskip
\section{$L^p$ convergence for the KFP model}
\label{sec:Conclude}
\setcounter{equation}{0}
\setcounter{theo}{0}

Before going to the proof of our main theorem, we need two last deduced results.
\begin{lem}\label{L61}
For any $\epsilon > 0$ small enough, we have
 \beqn
\nonumber
\Vert \AA S_\BB(t) \Vert_{L^2(G^{- (\frac 1 2 +\epsilon)}) \to L^2(G^{-(\frac 1 2 +\epsilon)})} \lesssim e^{- a t^{\frac {\gamma} {2-   \gamma }} } , \quad \forall t \ge 0,
\eeqn
and
\beqn
\nonumber
\Vert \AA \SS_B(t)  \Vert_{ L^1(G^{-(\frac 1 2 +\epsilon)})   \to L^1(G^{-(\frac 1 2 +\epsilon)})   }\lesssim e^{-a t^{\frac {\gamma} {2-   \gamma }}}, \quad \forall t \ge 0,
\eeqn
for some $a>0$. Similarly for any $0 < b <\frac {\gamma} {2-\gamma}$ and for any $\epsilon > 0$ small enough, we have
\beqn
\nonumber
\Vert \AA S_\BB(t) \Vert_{ L^1(G^{-(\frac 1 2 +\epsilon)})   \to L^2(G^{-(\frac 1 2 +\epsilon)}) } \lesssim t^{-\alpha}e^{- a t^{b}}, \quad \forall t \ge 0,
\eeqn
for $\alpha =\frac {5d+1} {2}$ and some $a > 0$.
\end{lem}

\noindent{\sl Proof of  Lemma \ref{L61}.}
 The first two inequalities are obtained obviously by Lemma \ref{L31} and the property of $\AA= M\chi_R$. For the third inequality we split it into two parts, $t \in [0, \eta]$ and $t > \eta$, where $\eta$ is defined in Theorem \ref{T41}. When $t \in [0, \eta]$ , we have $e^{-at^{\frac {\gamma} {2-   \gamma }}} \ge e^{-a \eta^{\frac {\gamma} {2-   \gamma }}}$, by Theorem \ref{T41}, we have
\beqn
\nonumber
\Vert  \AA \SS_B(t)  \Vert_{   L^1(G^{-(\frac 1 2 +\epsilon)})    \to L^2(G^{-(\frac 1 2 +\epsilon)})  }\lesssim t^{-\alpha} \lesssim t^{-\alpha}e^{-a t^{\frac {\gamma} {2-   \gamma }}} ,\quad \forall t \in [0, \eta],
\eeqn
for some $a > 0$. When $ t \ge \eta$,  by Theorem \ref{T41}, we have
\beqn
\nonumber
\Vert  S_\BB(\eta) \Vert_{  L^2(G^{-(\frac 1 2 +\epsilon)})   \to L^2(G^{-(\frac 1 2 +\epsilon)}) }\lesssim  \eta^\alpha \lesssim 1,
\eeqn
and by Lemma \ref{L31}
\beqn
\nonumber
\Vert \SS_\BB(t - \eta)  \Vert_{L^2(G^{-(\frac 1 2 +\epsilon)}) \to L^2(G^{-\frac 1 2 }) }\lesssim e^{-a (t-\eta)^{\frac {\gamma} {2-   \gamma }}} \lesssim  e^{-a t^{\frac {\gamma} {2-   \gamma }}},
\eeqn
gathering the two inequalities, we have
\beqn
\nonumber
\Vert \AA \SS_B(t)  \Vert_{L^1(G^{-1/2 (1 + 2\epsilon )}) \to L^2(G^{-1/2( 1+ 2\epsilon)})}\lesssim e^{-at^{\frac {\gamma} {2-   \gamma }}} \lesssim t^{-\alpha}e^{-a t^b}, \quad \forall t > \eta,
\eeqn
for any $0 < b < \frac \gamma {2-\gamma}$, the proof is ended by combining the two cases above.
\qed

\begin{lem}\label{L62}
Similarly as Lemma \ref{L61}. For any $p \in (2,\infty)$, we have
\beqn
\nonumber
\Vert S_\BB(t)  \AA \Vert_{L^2(G^{-1/2}) \to L^2(G^{-1/2 })} \lesssim e^{- a t^{\frac {\gamma} {2-   \gamma }} } , \quad \forall t \ge 0.
\eeqn
and
\beqn
\nonumber
\Vert  \SS_B(t) \AA \Vert_{L^p(G^{-1/2 })\to L^p(G^{-1/2 
})}\lesssim e^{-a t^{\frac {\gamma} {2-   \gamma }}}, \quad \forall t \ge 0.
\eeqn
for some $a>0$. And for any $0 < b < \frac \gamma {2-\gamma}$ we have
\beqn
\nonumber
\Vert S_\BB(t) \AA \Vert_{L^2(G^{-1/2 }) \to L^p(G^{-1/2 }) } \lesssim t^{-\beta}e^{- a t^b}, \quad \forall t \ge 0.
\eeqn
for some $\beta >0$ and some $a > 0$.
\end{lem}

\noindent The proof of Lemma \ref{L62} is similar to the proof of Lemma \ref{L61} and thus skipped.

\begin{lem}\label{L63}
let $X, Y$ be two Banach spaces, $S(t)$ a semigroup such that for  all $t \ge 0 $and some $0 < a,0 < b<1$ we have
\beqn
\nonumber
\Vert S(t) \Vert_{X \to X} \le C_X e^{-at^b}, \ \ \Vert S(t) \Vert_{Y \to Y} \le C_Y e^{-at^b},
\eeqn
and for some $0 < \alpha $, we have 
\beqn
\nonumber
\Vert S(t) \Vert_{X \to Y} \le C_{X,Y} t^{-\alpha} e^{-at^b}.
\eeqn
Then we can have that for all integer $ n  > 0$
\beqn
\nonumber
\Vert S^{(*n)}(t) \Vert_{X \to X} \le C_{X,n} t^{n-1} e^{-at^b},
\eeqn
similarly
\beqn
\nonumber
\Vert S^{(*n)}(t) \Vert_{Y \to Y} \le C_{Y,n} t^{n-1} e^{-at^b},
\eeqn
and
\beqn
\nonumber
\Vert S^{(*n)}(t) \Vert_{X \to Y} \le C_{X,Y,n} t^{n-\alpha-1} e^{-at^b}.
\eeqn
In particular for $\alpha+1 < n$, and for any $ b^{*} < b$
\beqn
\nonumber
\Vert S^{(*n)}(t) \Vert_{X \to Y} \le C_{X,Y,n} e^{-at^{b^*}}.
\eeqn
\end{lem}

\noindent { \sl Proof of Lemma \ref{L63}.}
The proof is the same as  Lemma 2.5 in \cite{MQT}, plus the fact $t^b \le s^b +(t-s)^b$ for any $0 \le s \le t, 0 < b < 1$.\qed

Then we come to the final proof.

\noindent {Proof of Theorem \ref{T11}.}
We only prove the case when $m=G^{\frac {p-1} {p}(1+\epsilon)}, \quad p\in[1, 2]$, for the proof of the other cases, one need only replace the use of Lemma \ref{L61} in the following proof by Lemma \ref{L62} and Theorem \ref{T41}. We will prove $p=1$ first, this time we need to prove
\beqn
\nonumber
\Vert S_\LL(I - \Pi)(t) \Vert_{L^{1}(G^{-\epsilon}) \to L^1} \lesssim  e^{-at^{b}},
\eeqn
for any $0< b < \frac {\gamma} {2-\gamma}$, where $I$ is the identity operator and $\Pi$ is a projection operator defined by
\beqn
\nonumber
\Pi (f) =\MM(f) G.
\eeqn 
First, Iterating the Duhamel's formula we split it into 3 terms
\bear
\nonumber
S_\LL(I-\Pi) &=&(I-\Pi)\{S_\BB + \sum_{l=1}^{n-1}( S_\BB \AA)^{(*l)}* (S_\BB) \}
\\ \nonumber
&&+\{ (I-\Pi)S_\LL \}*(\AA S_\BB(t))^{*n},
\eear
and we will estimate them separately. By Lemma \ref{L31}, we have
\beqn
\nonumber
\Vert S_\BB(t) \Vert_{L^{1}(G^{-\epsilon}) \to L^1} \lesssim  e^{-at^{\frac {\gamma} {2-   \gamma }}},
\eeqn
the first term is thus estimated. For the second term, still using Lemma \ref{L31}, we get
\beqn
\nonumber
\Vert S_\BB(t)\AA \Vert_{L^1(G^{-\epsilon}) \to L^1} \lesssim e^{-at^{\frac {\gamma} {2-   \gamma }}},
\eeqn
by Lemma \ref{L63}, we have
\beqn
\nonumber
\Vert (S_\BB(t)\AA)^{*l} \Vert_{L^1(G^{-\epsilon}) \to L^1} \lesssim t^{l-1} e^{-at^{\frac {\gamma} {2-   \gamma }}},
\eeqn
thus the second term is estimated. For the last term  by Lemma \ref{L31}
\beqn
\nonumber
\Vert \AA S_\BB(t) \Vert_{L^1( G^{-\epsilon}) \to L^1(G^{-(\frac 1 2 + \epsilon)})} \lesssim e^{- a t^{\frac {\gamma} {2-   \gamma }}} .
\eeqn
By Lemma \ref{L61} and \ref{L63}, for any $0 < b < \frac {\gamma } {2-\gamma}$, we have
\beqn
\nonumber
\Vert (\AA S_\BB)^{(*n-1)}(t) \Vert_{L^1(G^{-(\frac 1 2 + \epsilon)}) \to L^2(G^{-(\frac 1 2 + \epsilon)})   } \lesssim t^{n-\alpha-2}e^{- a t^b} ,
\eeqn
finally by Theorem \ref{T31}, we have
\beqn
\nonumber
\Vert S_\LL(t)(I-\Pi) \Vert_{  L^2(G^{-(\frac 1 2 + \epsilon)})   \to L^2(G^{-1/2}) } \lesssim e^{-at^{\frac {\gamma} {2-   \gamma }}}.
\eeqn
Taking $n > \alpha+2 $  the third term is estimated thus the proof of case $p=1$ is concluded by gathering the inequalities above. As the case $p=2$ ia already proved in Theorem $\ref{T31}$, the case $p \in (1, 2)$ follows by interpolation. \qed

\end{document}